\documentclass[11pt]{article}
\usepackage{amssymb}
\usepackage{amsmath}
\usepackage{amsthm}
\newcommand{\cA}{{\cal A}}
\newcommand{\cB}{{\cal B}}

\newcommand{\cH}{{\cal H}}
\newcommand{\cE}{{\cal E}}

\newcommand{\cI}{{\cal I}}
\newcommand{\cJ}{{\cal J}}
\newcommand{\cO}{{\cal O}}
\newcommand{\cL}{{\cal L}}
\newcommand{\cM}{{\cal M}}

\newcommand{\cF}{{\cal F}}

\newcommand{\cS}{{\cal S}}

\newcommand{\cW}{{\cal W}}
\newcommand{\cX}{{\cal X}}
\newcommand{\cY}{{\cal Y}}

\newcommand{\ZZ}{{\mathbb Z}}

\newcommand{\on}{\operatorname}

\newcommand{\Qlb}{\mathbb{\bar Q}_\ell}

\newcommand{\A}{\mathbb{A}}

\newcommand{\toup}[1]{\stackrel{#1}{\to}}
\newcommand{\hook}[1]{\stackrel{#1}{\hookrightarrow}}

\newcommand{\getsup}[1]{\stackrel{#1}{\gets}}
\newcommand{\Sp}{\on{\mathbb{S}p}}

\newcommand{\Sym}{\on{Sym}}
\newcommand{\SO}{\on{S\mathbb{O}}}
\newcommand{\Ker}{\on{Ker}}

\newcommand{\Aut}{\on{Aut}}

\newcommand{\RG}{\on{R\Gamma}}

\newcommand{\Bun}{\on{Bun}}

\newcommand{\Bunt}{\on{\widetilde\Bun}}

\newcommand{\Spec}{\on{Spec}}

\newcommand{\Gr}{\on{Gr}}

\newcommand{\GL}{\on{GL}}

\newcommand{\Fr}{{\on{Fr}}}

\newcommand{\pr}{\on{pr}}
\newcommand{\id}{\on{id}}

\newcommand{\QED}{$\square$} 
\newcommand{\Fq}{\mathbb{F}_q}  
\newcommand{\Fp}{\mathbb{F}_p}  
\newcommand{\iso}{{\widetilde\to}}

\newcommand{\comp}{\circ}
\newcommand{\Four}{\on{Four}}
\renewcommand{\H}{{\on{H}}}   

\newcommand{\D}{\on{D}}       
\newcommand{\wt}{\widetilde}

\newcommand{\select}[1]{{\it{#1}}}

\renewcommand{\P}{{\on{P}}}

\newcommand{\<}{\langle}
\renewcommand{\>}{\rangle}

\newcommand{\Res}{\on{Res}}

\newcommand{\act}{\on{act}}
\newcommand{\dimrel}{\on{dim.rel}}

\renewcommand{\Im}{\on{Im}}

\newcommand{\twolim}{\mathop{2\!\!-\!\!{\rm{lim}}}}
\newcommand{\diag}{\mathop{\vartriangle}}
\newcommand{\mult}{\on{mult}}
\newcommand{\inv}{\on{inv}}
\newcommand{\isoup}[1]{\stackrel{#1}{\iso}}
\newcommand{\Shr}{\on{Shr}}

\newtheorem{Lm}{Lemma}
\newtheorem{Th}{Theorem}
\newtheorem{Pp}{Proposition}

\newtheorem{Slm}{Sublemma}

\theoremstyle{remark}
\newtheorem{Rem}{Remark}

\theoremstyle{definition}
\newtheorem{Def}{Definition}

\newenvironment{Prf}{\par\noindent {\it Proof }}{\QED}
\newcommand{\Step}[1]{\par\noindent{\bf Step {#1}}.}

\begin{document}
\author{Vincent Lafforgue, Sergey Lysenko}
\title{Geometric Weil representation: local field case}
\date{}
\maketitle
\begin{abstract}
\noindent{\scshape Abstract}\hskip 0.8 em  
Let $k$ be an algebraically closed field of characteristic $>2$, $F=k((t))$ and $G=\Sp_{2d}$. In this paper we propose a geometric analog of the Weil representation of the metaplectic group $\wt G(F)$. This is a category of certain perverse sheaves on some stack, on which  $\wt G(F)$ acts by functors. This construction will be used in \cite{L2} (and subsequent publications) for the proof of the geometric Langlands functoriality for some dual reductive pairs.   
\end{abstract}

\bigskip\bigskip

{\centerline{\scshape 1. Introduction}}

\bigskip\noindent
1.1 This paper followed by \cite{L2} form a series, where we prove the geometric Langlands functoriality for the dual reductive pair $\Sp_{2n}, \SO_{2m}$ (in the everywhere nonramified case). 

 Let $k=\Fq$ with $q$ odd, set $\cO=k[[t]]\subset F=k((t))$. Write $\Omega$ for the completed module of relative differentials of $\cO$ over $k$. Let $M$ be a free $\cO$-module of rank $2d$ with symplectic form $\wedge^2 M\to\Omega$, set $G=\Sp(M)$. The group $G(F)$ admits a nontrivial metaplectic extension 
$$
1\to \{\pm 1\}\to\tilde G(F)\to G(F)\to 1 
$$
(defined up to a unique isomorphism). Let $\psi: k\to\Qlb^*$ be a nontrivial additive character, let $\chi: \Omega(F)\to\Qlb^*$ be given by $\chi(\omega)=\psi(\Res\omega)$. Write $H=M\oplus\Omega$ for the Heisenberg group of $M$ with operation 
$$
(m_1, a_1)(m_2, a_2)=(m_1+m_2, \, a_1+a_2+\frac{1}{2}\omega\<m_1,m_2\>)\;\;\;\;\; m_i\in M, a_i\in\Omega
$$ 
Denote by $\cS_{\psi}$ the Weil representation of $H(M)(F)$ with central character $\chi$. As a representation of $\tilde G(F)$, it decomposes $\cS_{\psi}\,\iso\,\cS_{\psi, odd}\oplus \cS_{\psi, even}$ into a direct sum of two irreducible smooth representations, where the even (resp., the odd) part is unramified (resp., ramified). 

 The discovery of this representation by A.~Weil in \cite{W} had a major influence on the theory of automorphic forms (among numerous developpements and applications are Howe duality for reductive dual pairs, particular cases of classical Langlands functoriality, Siegel-Weil formulas, relation with L-functions, representation-theoretic approach to the theory of theta-series.  We refer the reader to \cite{Ge}, \cite{LV}, \cite{H}, \cite{MVW}, \cite{Pr} for history and details). 

 In this paper we introduce a geometric analog of the Weil representation $\cS_{\psi}$. The pioneering work in this direction is due to P. Deligne \cite{D}, where a geometric approach to the Weil representation of a symplectic group over a finite field was set up. It was further extended by Gurevich-Hadani in \cite{GH, GH1}. 
The point of this paper is to develop the geometric theory in the case when a finite field is replaced by a local non-archimedian field. 
  
 First, we introduce a $k$-scheme $\cL_d(M(F))$ of discrete lagrangian lattices in $M(F)$ and a certain $\mu_2$-gerb $\wt\cL_d(M(F))$ over it. We view the metaplectic group $\tilde G(F)$ as a group stack over $k$. We construct a category 
$$
W(\wt\cL_d(M(F)))
$$ 
of certain perverse sheaves on $\wt\cL_d(M(F))$, which provides a geometric analog of $\cS_{\psi,even}$. The metaplectic group $\tilde G(F)$ acts on the category  $W(\wt\cL_d(M(F)))$ by functors. This action is \select{geometric} in the sense that it comes from a natural action of $\tilde G(F)$ on $\wt\cL_d(M(F))$ (cf. Theorem~\ref{Th_2}). 
 
  The category $W(\wt\cL_d(M(F)))$ has a distinguished object $S_{M(F)}$ corresponding to the unique non-ramified vector of $\cS_{\psi, even}$. 
 
  Our category $W(\wt\cL_d(M(F)))$ is obtained from Weil representations of symplectic groups $\Sp_{2r}(k)$ by some limit procedure. This uses a construction of geometric canonical interwining operators for such representations. A similar result has been announced by Gurevich and Hadani in \cite{GH} and proved for $d=1$ in \cite{GH1}. We give a proof for any $d$ (cf. Theorem~\ref{Th_appendix_CIO}). When this paper has already been written we learned about a new preprint \cite{GH2}, where a result similar to our Theorem~1 is claimed to be proved for all $d$. However, the sheaves of canonical interwining operators constructed in \select{loc.cit.} and in this paper live on different bases.
    
  Finally, in Section~7 we give a global application. 
Let $X$ be a smooth projective curve. Write $\Omega_X$ for the canonical line bundle on $X$. Let
$G$ denote the sheaf of automorphisms of $\cO_X^d\oplus \Omega^d_X$ preserving the natural symplectic form 
$\wedge^2(\cO_X^d\oplus \Omega^d_X)\to\Omega_X$. 

 Our Theorem~\ref{Th_3} relates $S_{M(F)}$ with the theta-sheaf $\Aut$ on the moduli stack $\Bunt_G$ of metaplectic bundles on $X$ introduced in \cite{L1}. This result will play an important role in \cite{L2}.

\bigskip\noindent
1.2  {\scshape Notation\ } In Section~2 we let $k=\Fq$ of characteristic $p>2$. Starting from Section~3 we assume $k$ either finite as above or algebraically closed with a fixed inclusion $\Fq\hook{} k$. All the schemes (or stacks) we consider are defined over $k$. 

  Fix a prime $\ell\ne p$. For a scheme (or stack) $S$ write $\D(S)$ for the bounded derived category of $\ell$-adic \'etale sheaves on $S$, and $\P(S)\subset \D(S)$ for the category of perverse sheaves. 
  
 Fix a nontrivial character $\psi: \Fp\to\Qlb^*$, write
$\cL_{\psi}$ for the corresponding Artin-Shreier sheaf on $\A^1$. Fix a square root $\Qlb(\frac{1}{2})$ of the sheaf $\Qlb(1)$ on $\Spec \Fq$. Isomorphism classes of such correspond to square roots of $q$ in $\Qlb$. 

 If $V\to S$ and $V^*\to S$ are dual rank $n$ vector bundles over a stack $S$, we normalize 
the Fourier transform $\Four_{\psi}: \D(V)\to\D(V^*)$ by 
$\Four_{\psi}(K)=(p_{V^*})_!(\xi^*\cL_{\psi}\otimes p_V^*K)[n](\frac{n}{2})$,  
where $p_V, p_{V^*}$ are the projections, and $\xi: V\times_S V^*\to \A^1$ is the pairing.

 Our conventions about $\ZZ/2\ZZ$-gradings are those of \cite{L1}. 

\bigskip\bigskip

\newpage

\centerline{\scshape 2. Canonical interwining operators: finite field case}

\bigskip\noindent
2.1 Let $M$ be a symplectic $k$-vector space of dimension $2d$. The symplectic form on $M$ is denoted $\omega\<\cdot,\cdot\>$. 
The Heisenberg group $H=M\times \A^1$ with operation
$$
(m_1, a_1)(m_2, a_2)=(m_1+m_2, \, a_1+a_2+\frac{1}{2}\omega\<m_1,m_2\>)\;\;\;\;\; m_i\in M, a_i\in\A^1
$$ 
is algebraic over $k$. Set $G=\Sp(M)$. Write $\cL(M)$ for the variety of lagrangian subspaces in $M$.  Fix a 
one-dimensional $k$-vector space $\cJ$ (purely of degree $d\!\!\mod 2$ as $\ZZ/2\ZZ$-graded). Let $\cA$ be the (purely of degree zero as $\ZZ/2\ZZ$-graded) line bundle over $\cL(M)$ with fibre $\cJ\otimes\det L$ at $L\in\cL(M)$. Write $\wt\cL(M)$ for the gerb of square roots of $\cA$. The line bundle $\cA$ is $G$-equivariant, so $G$ acts naturally on $\wt\cL(M)$. 
 
 For a $k$-point $L\in\cL(M)$ write $L^0$ for a $k$-point of $\wt\cL(M)$ over $L$. Write 
$$
\bar L=L\oplus k,
$$
this is a subgroup of $H(k)$ equipped with the character $\chi_L: \bar L\to\Qlb^*$ given by $\chi_L(l,a)=\psi(a)$, $l\in L, a\in k$. Write 
$$
\cH_L=\{f:H(k)\to\Qlb\mid f(\bar l h)=\chi_L(\bar l)f(h),\;\;\mbox{for}\;\; \bar l\in \bar L, h\in H\}
$$
This is a representation of $H(k)$ by right translations. 
Write $\cS(H)$ for the space of all $\Qlb$-valued functions on $H(k)$. The group $G$ acts naturally in $\cS(H)$. For $L\in\cL(M), g\in G$ we have an isomorphism $\cH_L\to \cH_{gL}$ sending $f$ to $gf$. 

 The purpose of Sections~2 and 3 is to study the canonical interwining operators (and their geometric analogs) between various models $\cH_L$ of the Weil representation. The corresponding results for a finite field were formulated by Gurevich and Hadani \cite{GH} without a proof (we give all proofs for the sake of completeness). Besides, our setting is a bit different from \select{loc.cit}, we work with gerbs instead of the total space of the corresponding line bundles. 

\medskip\noindent
2.2 For $k$-points $L^0,N^0\in\wt\cL(M)$ we will define a canonical interwining operator 
$$
F_{N^0,L^0}: \cH_L\to \cH_N
$$ 
They will satisfy the properties
\begin{itemize}
\item $F_{L^0, L^0}=\id$
\item $F_{R^0,N^0}\comp F_{N^0,L^0}=F_{R^0,L^0}$ for any $R^0,N^0,L^0\in\wt\cL(M)$
\item for any $g\in G$ we have $g\comp F_{N^0,L^0}\comp g^{-1}=F_{gN^0, gL^0}$.
\item under the natural action of $\mu_2$ on the set $\wt\cL(M)(k)$ of (isomorphism classes of) $k$-points, 
$F_{N^0,L^0}$ is odd as a function of $N^0$ and of $L^0$.
\end{itemize}

In (Remark~\ref{Rem_function_F_classical}, Section~3.1) we will define a function $F^{cl}$ on the set of $k$-points of $\wt\cL(M)\times \wt\cL(M)\times H$, which we denote $F_{N^0,L^0}(h)$ for $h\in H$.  
It will realize the operator $F_{N^0,L^0}$ by
$$
(F_{N^0,L^0}f)(h_1)=\int_{h_2\in H} F_{N^0,L^0}(h_1h_2^{-1})f(h_2)dh_2
$$
All our measures on finite sets are normalized by requiring the volume of a point to be one. Given two functions $f_1,f_2: H\to\Qlb$ their convolution $f_1\ast f_2: H\to\Qlb$ is defined by 
$$
(f_1\ast f_2)(h)=\int_{v\in H} f_1(hv^{-1})f_2(v)dv\;\;\;\; h\in H
$$

  The function $F_{N^0,L^0}$ will satisfy the following:
 
\smallskip 
 
\begin{itemize}
\item $F_{N^0,L^0}(\bar n h\bar l)=\chi_N(\bar n)\chi_L(\bar l) F_{N^0,L^0}(h)\;\;$ for $\bar l\in\bar L, \bar n\in\bar N, h\in H$.

\item $F_{gN^0, gL^0}(gh)=F_{N^0, L^0}(h)$ for $g\in G$, $h\in H$.

\item Convolution property:
$F_{R^0, L^0}=F_{R^0, N^0}\ast F_{N^0, L^0}\;$  for any $R^0,N^0,L^0\in\wt\cL(M)$. 
\end{itemize}

\smallskip\noindent
2.3 First, we define the non-normalized function $\tilde F_{N, L}: H\to \Qlb$, it will depend only on $N,L\in\cL(M)$, not ot their inhanced structure. 

 Given $N,L\in\cL(M)$ let $\chi_{NL}: \bar N \bar L\to \Qlb$ be the function given by 
$$
\chi_{NL}(\bar n \bar l)=\chi_N(\bar n)\chi_L(\bar l),
$$
it is correctly defined.  Note that $\bar N\bar L=\bar L\bar N$ but $\chi_{NL}\ne \chi_{LN}$ in general. Set 
$$
\tilde F_{N, L}(h)=\left\{ 
\begin{array}{cl}
\chi_{NL}(h),   & \mbox{if}\;\; h\in \bar N\bar L\\
0, & \mbox{otherwise}
\end{array}   
\right.
$$ 
Note that $\chi_{LL}=\chi_L$. 

 Given $L,R,N\in\cL(M)$ with $N\cap L=N\cap R=0$, define $\theta(R,N,L)\in\Qlb$ as follows. There is a unique map $b: L\to N$ such that $R=\{l+b(l)\in L\oplus N\mid l\in L\}$. Set
$$
\theta(R,N,L)=\int_{l\in L} \psi(\frac{1}{2}\omega\<l, b(l)\>)dl
$$
This expression has been considered in (\cite{L1}, Appendix~B).

\begin{Lm} 
\label{Lm_appendix_CIT}
1) Let $L,N\in\cL(M)$. If $L\cap N=0$ then $\tilde F_{L,N}\ast \tilde F_{N,L}=q^{2d+1}\tilde F_{L,L}$. 

\medskip\noindent
2) Let $L,R,N\in\cL(M)$ with $N\cap L=N\cap R=0$. Then $\tilde F_{R,N}\ast \tilde F_{N,L}=q^{d+1}\theta(R,N,L) \tilde F_{R,L}$
\end{Lm}
\begin{Prf}
2) Using $L\oplus N=N\oplus R=M$, for $h\in H$ we get 
$$
(\tilde F_{R,N}\ast \tilde F_{N,L})(h)=q^{d+1}\int_{v\in \bar N\backslash H} \chi_{RN}(hv^{-1})\chi_{NL}(v)dv=
q^{d+1}\int_{r\in R} \chi_{RN}(h(-r,0))\chi_{NL}(r,0)dr
$$
Because of the equivariance property of $\tilde F_{R,N}\ast \tilde F_{N,L}$, we may assume $h=(n,0), n\in N$. We get 
\begin{multline}
\label{convolution_expression_one_appendix}
(\tilde F_{R,N}\ast \tilde F_{N,L})(h)=q^{d+1} \int_{r\in R} \chi_{RN}((n,0)(-r,0))\chi_{NL}(r,0)dr\\
=q^{d+1} \int_{r\in R} 
\psi(\omega\<r,n\>)\chi_{NL}(r,0)dr
\end{multline}

 The latter formula essentially says that the resulting function on $N$ is the Fourier transform of some local system on $R$ (the symplectic form on $M$ induces an isomorphism $R\,\iso\, N^*$). This will be used for geometrization in Lemma~\ref{Lm_two_appendix}. 
 
 There is a unique map $b: L\to N$ such that $R=\{l+b(l)\in L\oplus N\mid l\in L\}$. So, the above integral rewrites
\begin{multline}
\label{value_h_1_h_2_for_Lm}
(\tilde F_{R,N}\ast \tilde F_{N,L})(h)=q^{d+1}\int_{l\in L}
\psi(\omega\<l,n\>)\chi_{NL}((l+b(l),0)dl=\\
q^{d+1}\int_{l\in L}
\psi(\omega\<l,n\>)\chi_{NL}((b(l),\frac{1}{2}\omega\<l, b(l)\>)(l,0))dl=q^{d+1}\int_{l\in L}
\psi(\omega\<l,n\>+\frac{1}{2}\omega\<l, b(l)\>)dl
\end{multline}

Note that if $R=L$ then $b=0$ and the latter formula yields 1). 

Let us identify $N\,\iso L^*$ via the map sending $n\in N$ to the linear functional $l\mapsto \omega\<l,n\>$. Denote by $\<\cdot,\cdot\>$ the symmetric pairing between $L$ and $L^*$. By Sublemma~\ref{Slm_1} below, the value (\ref{value_h_1_h_2_for_Lm}) vanishes unless $n\in (R+L)\cap N=\Im b$. In the latter case pick $l_1\in L$ with $b(l_1)=n$. Then 
$$
\chi_{RL}(n,0)=\psi(-\frac{1}{2}\omega\<l_1, b(l_1)\>)
$$ 
So, we get for $L'=\Ker b$
$$
(\tilde F_{R,N}\ast \tilde F_{N,L})(h)=q^{d+1+\dim L'}\chi_{RL}(h) \int_{l\in L/L'} \psi(\frac{1}{2}\omega\<l, b(l)\>)dl
$$
We are done.
\end{Prf}

\begin{Slm} 
\label{Slm_1}
Let $L$ be a $d$-dimensional $k$-vector space, $b\in\Sym^2 L^*$ and $u\in L^*$. View $b$ as a map $b: L\to L^*$, let $L'$ be the kernel of $b$. Then 
\begin{equation}
\label{expression_gauss_int_linear_term}
\int_{l\in L} \psi(\<l,u\>+\frac{1}{2}\<l, b(l)\>)dl
\end{equation}
is supported at $u\in (L/L')^*$ and there equals
$$
q^{\dim L'} \psi(-\frac{1}{2}\<b^{-1}u, u\>) 
\int_{L/L'}\psi(\frac{1}{2}\<l, b(l)\>)dl, 
$$
where $b: L/L'\,\iso\, (L/L')^*$, so that $b^{-1}u \in L/L'$. 
(Here the scalar product is between $L$ and $L^*$, so is symmetric).
\end{Slm}
\begin{Prf}
Let $L'\subset L$ denote the kernel of $b:L\to L^*$. Integrating first along the fibres of the projection  
$L\to L/L'$ we will get zero unless $u\in (L/L')^*$.
For any $l_0\in L$ the integral (\ref{expression_gauss_int_linear_term})
equals
\begin{multline*}
\int_{l\in L}\psi(\<l+l_0, u\>+\frac{1}{2}\<l+l_0, b(l)+b(l_0)\>)dl=\psi(\<l_0,u\>+\frac{1}{2}\<l_0, b(l_0)\>)\int_{l\in L} \psi(\<l,u+b(l_0)\>+\frac{1}{2}\<l, b(l)\>)dl
\end{multline*}
Assuming $u\in (L/L')^*$ take $l_0$ such that $u=-b(l_0)$. Then (\ref{expression_gauss_int_linear_term}) becomes
$$
\psi(\frac{1}{2}\<l_0, u\>)\int_{l\in L} \psi(\frac{1}{2}\<l, b(l)\>)dl
$$
We are done.
\end{Prf}

\medskip

\begin{Rem} The expression (\ref{expression_gauss_int_linear_term}) is the Fourier transform from $L$ to $L^*$. In the geometric setting we will use 2) of Lemma~\ref{Lm_appendix_CIT} only under the additional assumption $R\cap L=0$. 
\end{Rem}

\bigskip\bigskip

\centerline{\scshape 3. Geometrization}

\bigskip\noindent
3.1 Let $M$, $H$, $\cL(M)$ and $\wt\cL(M)$ be as in Section~2.1. Remind that $G=\Sp(M)$. For each $L\in \cL(M)$ we have a rank one local system $\chi_L$ on $\bar L=L\times \A^1$ defined by $\chi_L=\pr^*\cL_{\psi}$, where $\pr: L\times \A^1\to \A^1$ is the projection. 
Let $\cH_L$ denote the category of perverse sheaves on $H$ which are $(\bar L, \chi_L)$-equivariant under the left multiplication, this is a full subcategory in $\P(H)$. Write $\D\!\cH_L\subset \D(H)$ for the full subcategory of objects whose all perverse cohomologies lie in $\cH_L$. 

 Denote by $C\to \cL(M)$ (resp., $\bar C\to \cL(M)$) the vector bundle whose fibre over $L\in\cL(M)$ is $L$ (resp., $\bar L=L\times\A^1$). Its inverse image to $\wt\cL(M)$ is denoted by the same symbol. 
 
 Write $\chi_{\bar C}$ for the local system $p^*\cL_{\psi}$ on $\bar C$, where $p: \bar C\to\A^1$ is the projection on the center sending $(L\in \cL(M), (l,a)\in\bar L)$ to $a$. Consider the maps
$$
\pr, \act_{lr}: \bar C\times \bar C\times H\to \cL(M)\times \cL(M)\times H\times H
$$ 
where $\act_{lr}$ sends $(\bar n\in \bar N, \bar l\in\bar L, h)$ to $(N,L, \bar n h \bar l)$, and $\pr$ sends the above point to $(N,L, h)$. We say that a perverse sheaf $K$ on $\cL(M)\times \cL(M)\times H$ is \select{$\act_{lr}$-equivariant} if it admits an isomorphism
$$
\act_{lr}^*K\,\iso\, \pr^*K\otimes \pr_1^*\chi_{\bar C}\otimes\pr_2^*\chi_{\bar C}
$$
satisfying the usual associativity condition and whose restriction to the unit secton is the identity (such isomorphism is unique if it exists). One has a similar definition for  $\wt\cL(M)\times \wt\cL(M)\times H$. 
 
 Let 
$$
\act_G: G\times \wt\cL(M)\times \wt\cL(M)\times H\to \wt\cL(M)\times \wt\cL(M)\times H
$$
be the action map sending $(g, N^0, L^0, h)$ to 
$$
(gN^0, gL^0, gh)
$$ 
For this map we have a usual notion of a \select{$G$-equivariant perverse sheaf} on $\wt\cL(M)\times \wt\cL(M)\times H$. As $G$ is connected, a perverse sheaf on $\wt\cL(M)\times \wt\cL(M)\times H$ admits at most one $G$-equivariant structure. 

 If $S$ is a stack then for $K,F\in \D(S\times H)$ define their convolution $K\ast F\in \D(S\times H)$ by
$$
K\ast F=\mult_!(\pr_1^*K\otimes \pr_2^*F)\otimes(\Qlb[1](\frac{1}{2}))^{d+1-2\dim\cL(M)},   
$$
here $\pr_i:  S\times H\times H\to S\times H$ is the projection to the $i$-th component in the pair $H\times H$ (and the identity on $S$). The multiplication map $\mult: H\times H\to H$ sends $(h_1,h_2)$ to $h_1h_2$. 

 Let 
\begin{equation}
\label{map_diag_one} 
(\cL(M)\times H)_{\diag}\hook{} \cL(M)\times H
\end{equation}
be the closed subscheme of those $(L\in\cL(M),h\in H)$ for which $h\in \bar L$. Let 
$$
\alpha_{\diag}:  
(\cL(M)\times H)_{\diag}\to\A^1
$$ 
be the map sending $(L,h)$ to $a$, where $h=(l,a)$, $l\in L, a\in\A^1$. Define a perverse sheaf
$$
\tilde F_{\diag}=\alpha_{\diag}^*\cL_{\psi}\otimes(\Qlb[1](\frac{1}{2}))^{d+1+\dim\cL(M)},
$$
which we extend by zero under (\ref{map_diag_one}). 

 Since $\wt\cL(M)\to\cL(M)$ is a $\mu_2$-gerb, $\mu_2$ acts on each $K\in \D(\wt\cL(M))$, and we say that $K$ is \select{genuine} if $-1\in\mu_2$ acts on $K$ as $-1$. 

\begin{Th} 
\label{Th_appendix_CIO}
There exists an irreducible perverse sheaf $F$ on $\wt\cL(M)\times\wt\cL(M)\times H$ (pure of weight zero) with the following properties:
\begin{itemize}
\item for the diagonal map $i: \wt\cL(M)\times H\to \wt\cL(M)\times\wt\cL(M)\times H$ the complex $i^*F$ identifies canonically with the inverse image of 
$$
\tilde F_{\diag}\otimes (\Qlb[1](\frac{1}{2}))^{\dim\cL(M)}
$$ 
under the projection $\wt\cL(M)\times\H\to \cL(M)\times H$. 
\item $F$ is $\act_{lr}$-equivariant;
\item $F$ is $G$-equivariant;
\item $F$ is genuine in the first and the second variable;
\item convolution property for $F$ holds, namely
for the $ij$-th projections 
$$
q_{ij}: \wt\cL(M)\times\wt\cL(M)\times\wt\cL(M)\times H\to \wt\cL(M)\times\wt\cL(M)\times H
$$
inside the triple $\wt\cL(M)\times\wt\cL(M)\times\wt\cL(M)$ we have 
$(q_{12}^*F)\ast (q_{23}^*F)\,\iso\, q_{13}^*F$ canonically. 
\end{itemize}
\end{Th}

\noindent
 The proof of Theorem~\ref{Th_appendix_CIO} is given in Sections~3.2-3.4.
 
\begin{Rem}
\label{Rem_function_F_classical}
In the case $k=\Fq$ define $F^{cl}$ as the trace of the geometric Frobenius on $F$.
\end{Rem}
 
\medskip\noindent
3.2 Let $U\subset \cL(M)\times \cL(M)$ be the open subset of pairs $(N,L)\in \cL(M)\times \cL(M)$ such that $N\cap L=0$. Define a perverse sheaf $\tilde F_U$ on $U\times H$ as follows. Let 
$$
\alpha_U: U\times H\to\A^1
$$ 
be the map sending $(N,L,h)$ to $a+\frac{1}{2}\omega\<l,n\>$, where $l\in L, n\in N, a\in\A^1$ are uniquely defined by $h=(n+l,a)$. Set
\begin{equation}
\label{def_tilde_F_U}
\tilde F_U=\alpha_U^*\cL_{\psi}\otimes(\Qlb[1](\frac{1}{2}))^{\dim H+2\dim\cL(M)}
\end{equation}
 
Write $U\times_{\cL(M)} U\subset \cL(M)\times\cL(M)\times \cL(M)$ for the open subscheme classifying $(R,N,L)$ with $N\cap L=N\cap R=0$. Let 
$$
q_i: U\times_{\cL(M)} U\to U
$$ 
be the projection on the $i$-th factor, so $q_1$ (resp., $q_2$) sends $(R,N,L)$ to $(R,N)$ (resp., to $(N,L$)). 
Let $q: U\times_{\cL(M)} U\to \cL(M)\times\cL(M)$ be the map sending $(R,N,L)$ to $(R,L)$. Write 
$$
(U\times_{\cL(M)} U)_0=q^{-1}(U)
$$ 
 
  The geometric analog of $\theta(R,N,L)$ is the following (shifted) perverse sheaf $\Theta$ on $U\times_{\cL(M)} U$. Let $\pi_C: C_3\to U\times_{\cL(M)} U$ be the vector bundle whose fibre over $(R,N,L)$ is $L$. We have a map $\beta: C_3\to\A^1$ defined as follows. Given a point $(R,N,L)\in U\times_{\cL(M)} U$, there is a unique map $b: L\to N$ such that $R=\{l+b(l)\in L\oplus N=M\mid l\in L\}$. Set $\beta(R,N,L, l)=\frac{1}{2}\omega\<l, b(l)\>$. Set
$$
\Theta=(\pi_C)_!\beta^*\cL_{\psi}\otimes(\Qlb[1](\frac{1}{2}))^d
$$ 

 Write $Y=\cL(M)\times\cL(M)$, let $\cA_Y$ be the ($\ZZ/2\ZZ$-graded purely of degree zero) line bundle on $Y$ whose fibre at $(R,L)$ is $\det R\otimes\det L$. Write
$\tilde Y$ for the gerb of square roots of $\cA_Y$. Note that $\cA_Y$ is $G$-equivariant, so $G$ acts on $\tilde Y$ naturally. 

The following perverse sheaf $S_M$ on $\tilde Y$ was introduced in (\cite{L1}, Definition~2). Let $Y_i\subset Y$ be the locally closed subscheme given by $\dim(R\cap L)=i$ for $(R,L)\in Y_i$. The restriction of $\cA_Y$ to each $Y_i$ admits the following $G$-equivariant square root. For a point $(R,L)\in Y_i$ we have an isomorphism  
$L/(R\cap L)\,\iso\, (R/(R\cap L))^*$ sending $l$ to the functional $r\mapsto \omega\<r,l\>$. It induces a $\ZZ/2\ZZ$-graded isomorphism $\det R\otimes\det L\,\iso\, \det (R\cap L)^2$. 

So, for the restriction $\tilde Y_i$ of the gerb $\tilde Y\to Y$ to $Y_i$ we get a trivialization
\begin{equation}
\label{iso_triv_gerb_Y_i}
\tilde Y_i\,\iso\, Y_i\times B(\mu_2)
\end{equation}
Write $W$ for the nontrivial local system of rank one on $B(\mu_2)$ corresponding to the covering $\Spec k\to B(\mu_2)$. 

\begin{Def} 
\label{Def_S_M}
Let $S_{M,g}$ (resp., $S_{M,s}$) denote the intermediate extension of
$$
(\Qlb\boxtimes W)\otimes(\Qlb[1](\frac{1}{2}))^{\dim Y}
$$
from $\tilde Y_0$ to $\tilde Y$ (resp., of $(\Qlb\boxtimes W)\otimes(\Qlb[1](\frac{1}{2}))^{\dim Y-1}$ from $\tilde Y_1$ to $\tilde Y$). Set $S_M=S_{M,g}\oplus S_{M,s}$.
\end{Def}

 Let 
$$
\pi_Y: U\times_{\cL(M)} U\to \tilde Y
$$ 
be the map sending $(R,N,L)$ to 
$$
(R,L, \cB, \epsilon: \cB^2\,\iso\, \det R\otimes\det L),
$$
where $\cB=\det L$ and $\epsilon$ is the isomorphism induced by $\epsilon_0$. Here
$\epsilon_0: L\,\iso\, R$ is the isomorphism sending $l\in L$ to $l+b(l)\in R$. In other words, $\epsilon_0$ sends $l$ to the unique $r\in R$ such that $r=l\mod N\in M/N$. 
Write also $\tilde U=\tilde Y_0$. 

 Define $\cE\in\D(\Spec k)$ by 
$$
\cE=\RG_c(\A^1, \beta_0^*\cL_{\psi})\otimes \Qlb[1](\frac{1}{2}),
$$
where $\beta_0:\A^1\to\A^1$ sends $x$ to $x^2$. Then $\cE$ is a 1-dimensional vector space placed in cohomological degree zero. The geometric Frobenius $\Fr_{\Fq}$ acts on $\cE^2$ by $1$ if $-1\in (\Fq^*)^2$ and by $-1$ otherwise. A choice of $\sqrt{-1}\in k$ yields an isomorphism $\cE^2\,\iso\,\Qlb$, so $\cE^4\,\iso\, \Qlb$ canonically. 

As in (\cite{L1}, Proposition~5), one gets a canonical isomorphism
\begin{equation}
\label{iso_from_AnnENS}
\pi^*_Y (S_{M,g}\otimes \cE^d\oplus S_{M,s}\otimes \cE^{d-1})\,\iso\, \Theta\otimes(\Qlb[1](\frac{1}{2}))^{2\dim\cL(M)}
\end{equation}
Since $d\ge 1$, the restriction $\pi_Y: (U\times_{\cL(M)} U)_0\to \tilde U$ is smooth of relative dimension $\dim\cL(M)$, with geometrically connected fibres. It is convenient to introduce a rank one local system $\Theta_U$ on $\tilde U$ equipped with a canonical isomorphism
\begin{equation}
\label{iso_Theta_descent}
\Theta\,\iso\, \pi_Y^*\Theta_U
\end{equation}
over $(U\times_{\cL(M)} U)_0$. The local system $\Theta_U$ is defined up to a unique isomorphism. 

 Let $i_U: U\to U\times_{\cL(M)} U$ be the map sending $(L,N)$ to $(L,N,L)$. Let $p_1: U\to \cL(M)$ be the projection sending $(L,N)$ to $L$. 

\begin{Lm} 
\label{Lm_two_appendix}
1) The complex 
$$
(q_1^*\tilde F_U)\ast (q_2^*\tilde F_U)\otimes(\Qlb[1](\frac{1}{2}))^{\dim\cL(M)}
$$
 is an irreducible perverse sheaf on $U\times_{\cL(M)} U\times H$ pure of weight zero. We have canonically
$$
i_U^*((q_1^*\tilde F_U)\ast (q_2^*\tilde F_U))\,\iso\, p_1^*\tilde F_{\diag}\otimes(\Qlb[1](\frac{1}{2}))^{\dim\cL(M)}
$$ 
over $U\times H$.\\
2) There is a canonical isomorphism 
$$
(q_1^*\tilde F_U)\ast (q_2^*\tilde F_U)\,\iso\, 
q^*\tilde F_U\otimes \Theta
$$
over $(U\times_{\cL(M)} U)_0\times H$.
\end{Lm}
\begin{Prf}
1) Follows from the properties of the Fourier transform as in Lemma~\ref{Lm_appendix_CIT}, formula (\ref{convolution_expression_one_appendix}). 

\smallskip\noindent
2) The proof of Lemma~\ref{Lm_appendix_CIT} goes through in the geometric setting. Our additional assumption that $(R,N,L)\in (U\times_{\cL(M)} U)_0$ means that $b:L\to N$ is an isomorphism (it simplifies the argument a little). 
\end{Prf} 

\medskip

\begin{Rem} Let $i_{\diag}: \cL(M)\to \tilde Y$ be the map sending $L$ to $(L,L,\cB=\det L)$ equipped with the isomorphism $\id: \cB^2\,\iso\,\det L\otimes\det L$. The commutative diagram
\begin{equation}
\label{diag_for_diag}
\begin{array}{ccc}
U & \toup{i_U} & U\times_{\cL(M)} U\\
\downarrow\lefteqn{\scriptstyle p_1} && \downarrow\lefteqn{\scriptstyle \pi_Y}\\
\cL(M) & \toup{i_{\diag}} & \tilde Y
\end{array}
\end{equation}
together with (\ref{iso_from_AnnENS}) yield a canonical isomorphism
$$
i_{\diag}^*S_M\,\iso\, \left\{
\begin{array}{cc}
\cE^{-d}\otimes(\Qlb[1](\frac{1}{2}))^{2\dim\cL(M)-d}, & d\;\mbox{is even}\\\\
\cE^{1-d}\otimes(\Qlb[1](\frac{1}{2}))^{2\dim\cL(M)-d}, & d\;\mbox{is odd}
\end{array}
\right.
$$
\end{Rem}

\smallskip\noindent
3.3  Consider the following diagram
$$
\begin{array}{ccc}
\tilde U\; \getsup{\tilde q_1} & (U\times_{\cL(M)} U)_0 & \toup{\tilde q_2} \;\tilde U\\
& \downarrow\lefteqn{\scriptstyle \tilde q}\\
& \tilde U
\end{array}
$$ 
Here  $\tilde q$ is the restriction of $\pi_Y$, and the map $\tilde q_i$ is the lifting of $q_i$ defined as follows.
We set $\tilde q_1(R,N,L)=\tilde q(R, L,N)$ and $\tilde q_2(R,N,L)=\tilde q(N,R,L)$. 
 
  The following property is a geometric counterpart of the way the Maslov index of $(R,N,L)$ changes under permutations of three lagrangian subspaces. 

\begin{Lm} 
\label{Lm_Maslov_index}
1) For $i=1,2$ we have canonically over $(U\times_{\cL(M)} U)_0$
$$
\tilde q_i^*\Theta_U\otimes \tilde q^*\Theta_U\,\iso\,\Qlb
$$ 
2) We have  
$\Theta_U^2\,\iso\, \cE^{2d}$ canonically, so $\Theta_U^4\,\iso\, \Qlb$ canonically.
\end{Lm}
\begin{Prf}
1) The two isomorphisms are obtained similarly, we consider only the case $i=2$. For a point $(R,N,L)\in (U\times_{\cL(M)} U)_0$ we have isomorphisms $b: L\,\iso\, N$ and $b_0: L\,\iso\, R$ such that $R=\{l+b(l)\mid l\in L\}$ and $N=\{l+b_0(l)\mid l\in L\}$. Clearly, $b_0(-l)=l+b(l)$ for $l\in L$. Let $\beta_2: L\times L\to \A^1$ be the map sending $(l,l_0)$ to $\frac{1}{2}\omega\<l, b(l)\>+\frac{1}{2}\omega\<l, b_0(l)\>$. We must show that
$$
\RG_c(L\times L, \beta_2^*\cL_{\psi})\,\iso\, \Qlb[2d](d)
$$
The quadratic form $(l,l_0)\mapsto \omega\<l, b(l)\>-\omega\<l_0, b(l_0)\>$ is hyperbolic on $L\oplus L$. Consider the isotopic subspace $Q=\{(l,l)\in L\times L\mid l\in L\}$. Integrating first along the fibres of the projection $L\times L\to (L\times L)/Q$ and then over $(L\times L)/Q$, one gets the desired isomorphism.

\noindent
2) This follows from (\ref{iso_from_AnnENS}). 
\end{Prf}

\medskip

 Define a perverse sheaf $F_U$ on $\tilde U\times H$ by
$$
F_U=\pr_1^*\Theta_U\otimes \tilde F_U,
$$
it is understood that we take the inverse image of $\tilde F_U$ under the projection $\tilde U\times H\to U\times H$ is the above formula.  Let $F$ be the intermediate extension of $F_U$ under the open immersion $\tilde U\times H\subset \tilde Y\times H$. 

\begin{Rem} In the case $d=0$ we have $H=\A^1$ and $\tilde Y=B(\mu_2)$. In this case by definition $F=W\boxtimes \cL_{\psi}\otimes\Qlb[1](\frac{1}{2})$ over $\tilde Y\times H=B(\mu_2)\times\A^1$. 
\end{Rem}

\medskip

 Combining Lemma~\ref{Lm_Maslov_index} and 2) of Lemma~\ref{Lm_two_appendix}, we get the following.

\begin{Lm} 
\label{Lm_one_more_App}
We have canonically $(\tilde q_1^*F_U)\ast (\tilde q_2^*F_U)\,\iso\, \tilde q^*F_U\otimes \cE^{2d}$ over $(U\times_{\cL(M)} U)_0\times H$.
\end{Lm} 

We have  a map $\xi: \wt\cL(M)\times\wt\cL(M)\to \tilde Y$ sending $(\cB_1, N, \cB_1^2\,\iso\,\cJ\otimes\det N; \cB_2, L, \cB_2^2\,\iso\, \cJ\otimes\det L)$ to $(N,L,\cB)$, where $\cB=\cB_1\otimes\cB_2\otimes\cJ^{-1}$ is equipped with the natural isomorphism $\cB^2\,\iso\, \det N\otimes\det L$. 
The restriction of $F$ under 
$$
\xi\times\id: \wt\cL(M)\times\wt\cL(M)\times H\to \tilde Y\times H
$$ 
is also denoted by $F$. Clearly, $F$ is an irreducible perverse sheaf of weight zero. 

 Consider the cartesian square
$$
\begin{array}{ccc}
(U\times_{\cL(M)} U)_0\times H &\hook{} &  (U\times_{\cL(M)} U)\times H\\
 \downarrow\lefteqn{\scriptstyle \pi_Y\times\id} &&  \downarrow\lefteqn{\scriptstyle \pi_Y\times\id}\\
 \tilde U\times H & \hook{} & \tilde Y\times H
\end{array}
$$
This diagram together with Lemma~\ref{Lm_two_appendix} yield a canonical isomorphism over $(U\times_{\cL(M)} U)\times H$
\begin{equation}
\label{formula_explicit_F_appendix}
(\pi_Y\times\id)^*F\,\iso\, (q_1^*\tilde F_U)\ast (q_2^*\tilde F_U)
\end{equation}
by intermediate extension from $(U\times_{\cL(M)} U)_0\times H$. This gives an explicit formula for $F$.

 Consider the diagram 
$$
\begin{array}{ccc}
U\times H & \toup{i_U\times\id} & U\times_{\cL(M)} U\times H\\
\downarrow\lefteqn{\scriptstyle p_1\times\id} && \downarrow\lefteqn{\scriptstyle \pi_Y\times\id}\\
\cL(M)\times H & \toup{i_{\diag}\times\id} & \tilde Y\times H
\end{array}
$$
obtained from (\ref{diag_for_diag}) by multiplication with $H$. By Lemma~\ref{Lm_two_appendix} and (\ref{formula_explicit_F_appendix}), we get canonically 
$$
(p_1\times\id)^*(i_{\diag}\times\id)^*F\,\iso\, (p_1\times\id)^*\tilde F_{\diag}\otimes(\Qlb[1](\frac{1}{2}))^{\dim\cL(M)}
$$
Since $\tilde F_{\diag}$ is perverse and $p_1$ has connected fibres, this isomorphism descends to a uniquely defined isomorphism
$$
(i_{\diag}\times\id)^*F\,\iso\, \tilde F_{\diag}\otimes(\Qlb[1](\frac{1}{2}))^{\dim\cL(M)}
$$

 By construction, $F$ is $\act_{lr}$-equivariant and $G$-equivariant (this holds for $F_U$ and this property is preserved by the intermediate extension). 

\medskip\noindent
3.4 To finish the proof of Theorem~\ref{Th_appendix_CIO}, it remains to establish the convolution property of $F$. We actually prove it in the following form. 

 Write $\tilde Y\times_{\cL(M)}\tilde Y$ for the stack classifying $R,N,L\in\cL(M)$, one dimensional $k$-vector spaces $\cB_1,\cB_2$ equipped with isomorphisms $\cB_1^2\,\iso\, \det R\otimes\det N$ and $\cB_2^2\,\iso\,\det N\otimes\det L$. We have a diagram
 $$
\begin{array}{ccc}
\tilde Y \getsup{\tau_1} & \tilde Y\times_{\cL(M)}\tilde Y & \toup{\tau_2} \tilde Y\\
& \downarrow\lefteqn{\scriptstyle \tau}\\
& \tilde Y,
\end{array}
$$
where $\tau_1$ (resp., $\tau_2$) sends the above collection to $(R, N,\cB_1)\in\tilde Y$ (resp., $(N,L,\cB_2)\in\tilde Y$). The map $\tau$ sends the above collection to $(R,L,\cB)$, where $\cB=\cB_1\otimes\cB_2\otimes(\det N)^{-1}$ is equipped with $\cB^2\,\iso\, \det R\otimes\det L$. 

\begin{Pp} There is a canonical isomorphism over $(\tilde Y\times_{\cL(M)}\tilde Y)\times H$
\begin{equation}
\label{key_iso_convolution_appendix}
(\tau_1^*F)\ast(\tau_2^*F)\,\iso\,\tau^*F
\end{equation}
\end{Pp}
\begin{Prf}
\Step 1 Consider the diagram 
$$
\begin{array}{ccc}
(U\times_{\cL(M)} U)_0 & \toup{\tilde q_1\times\tilde q_2} & (\tilde U\times_{\cL(M)}\tilde U)_0\\
&\searrow\lefteqn{\scriptstyle \tilde q} & \downarrow\lefteqn{\scriptstyle\tau}\\
&& \tilde U
\end{array}
$$
It becomes 2-commutative over $\Spec \Fq(\sqrt{-1})$. More precisely, for $K\in \D(\tilde U)$ we have a canonical isomorphism functorial in $K$
$$
\tilde q^*K\otimes\cE^{2d}\,\iso\, (\tilde q_1\times\tilde q_2)^*\tau^*K
$$

 Indeed, let $(R,N,L)$ be a $k$-point of $(U\times_{\cL(M)} U)_0$, let $(R,N,L,\cB_1,\cB_2)$ be its image under $\tilde q_1\times\tilde q_2$. So, $\cB_1=\det N$ and $\pi_Y(R,L,N)=(R,N,\cB_1)$, $\cB_2=\det L$ and $\pi_Y(N,R,L)=(N,L,\cB_2)$. Write
$$
\tau(R,N,L,\cB_1,\cB_2)=(R,L,\cB, \delta: \cB^2\,\iso\, \det R\otimes\det L)
$$
Write $\tilde q(R,N,L)=(R,L, \cB, \delta_0: \cB^2\,\iso\, \det R\otimes\det L)$. It suffices to show that $\delta_0=(-1)^d \delta$.  
 
 Let $\epsilon_1: N\,\iso\, R$ be the isomorphism sending $n\in N$ to $r\in R$ such that $r=n\!\mod L$. Write $\epsilon_2: L\,\iso\, N$ for the isomorphism sending $l\in L$ to $n\in N$ such that $l=n\!\mod R$. Let $\epsilon_0: L\,\iso\, R$ be the isomorphism sending $l\in L$ to $r\in R$ such that $r=l\!\mod N$. We get two isomorphisms
$$
\id\otimes\det\epsilon_0,\;
\det \epsilon_1\otimes\det\epsilon_2: \;\; \det N\otimes\det L\,\iso\, \det R\otimes\det N
$$
We must show that $\id\otimes\det\epsilon_0=(-1)^d \det \epsilon_1\otimes\det\epsilon_2$. Pick a base $\{n_1,\ldots, n_d\}$ in $N$. Define $r_i\in R, l_i\in L$ by $n_i=r_i+l_i$. Then 
$$
\epsilon_1(n_i)=r_i,\;\;\; \epsilon_2(l_i)=n_i,\;\;\; \epsilon_0(l_i)=-r_i
$$
So, $\epsilon_0(l_1\wedge\ldots\wedge l_d)=(-1)^d r_1\wedge\ldots\wedge r_d$. On the other hand, $\det \epsilon_1\otimes\det\epsilon_2$ sends 
$$
(n_1\wedge\ldots\wedge n_d)\otimes (l_1\wedge\ldots\wedge l_d)
$$ 
to $(r_1\wedge\ldots\wedge r_d)\otimes (n_1\wedge\ldots\wedge n_d)$. 

\smallskip
\Step 2  
The isomorphism (\ref{iso_triv_gerb_Y_i}) for $i=0$ yields $(\tilde U\times_{\cL(M)}\tilde U)_0\,\iso\, (U\times_{\cL(M)} U)_0\times B(\mu_2)\times B(\mu_2)$. 
The corresponding 2-automorphisms $\mu_2\times\mu_2$ of $(\tilde Y\times_{\cL(M)}\tilde Y)$ act in the same way on both sides of (\ref{key_iso_convolution_appendix}). Now
from Step~1 it follows that the isomorphism of Lemma~\ref{Lm_one_more_App} descends under $\tilde q_1\times \tilde q_2$ to the desired isomorphism
(\ref{key_iso_convolution_appendix}) over $(\tilde U\times_{\cL(M)}\tilde U)_0\times H$. 

\smallskip
\Step 3 To finish the proof it suffices to show that $(\tau_1^*F)\ast(\tau_2^*F)$ is perverse, the intermediate extension under the open immersion
$$
(\tilde U\times_{\cL(M)}\tilde U)_0\times H\subset (\tilde Y\times_{\cL(M)}\tilde Y)\times H
$$
Let us first explain the idea informally, at the level of functions. In this step for $(N,R,\cB)\in\tilde Y$ we denote by 
$
F_{N,R,\cB}: H\to\Qlb
$ 
the function trace of Frobenius of the sheaf $F$. 

 Given $(R,N,\cB_1)\in\tilde Y$ and $(N,L,\cB_2)\in\tilde Y$ pick any $S,T\in\cL(M)$ such that $(R,S,N)\in U\times_{\cL(M)}U$, $(N,T,L)\in U\times_{\cL(M)}U$ and $S\cap T=S\cap L=0$. Assuming 
$$
(R,N,\cB_1)=\pi_Y(R,S,N)\;\;\; \mbox{and}\;\;\; (N,L,\cB_2)=\pi_Y(N,T,L),
$$
by (\ref{formula_explicit_F_appendix}) we get
\begin{multline*}
F_{R,N,\cB_1}\ast F_{N,L,\cB_2}=(\tilde F_{R,S}\ast \tilde F_{S,N})\ast (\tilde F_{N,T}\ast \tilde F_{T,L})=
q^{d+1}\theta(S,N,T)\tilde F_{R,S}\ast \tilde F_{S,T}\ast \tilde F_{T,L}\\ \\
=q^{2d+2}\theta(S,N,T)\theta(S,T,L)\tilde F_{R,S}\ast\tilde F_{S,L}=q^{2d+2}\theta(S,N,T)\theta(S,T,L) F_{R,L,\cB},
\end{multline*}
where $(R,L,\cB)=\pi_Y(R,S,L)$. Now we turn back to the geometric setting.

\smallskip
\Step 4 Consider the scheme $\cW$ classifying $(R,S,N)\in U\times_{\cL(M)} U$ and $(N,T,L)\in U\times_{\cL(M)} U$ such that $S\cap T=S\cap L=0$. Let 
$$
\kappa: \cW\to \tilde Y\times_{\cL(M)}\tilde Y
$$ 
be the map sending the above point to $(R,N,L,\cB_1,\cB_2)$, where $(R,N,\cB_1)=\pi_Y(R,S,N)$ and $(N,L,\cB_2)=\pi_Y(N,T,L)$. The map $\kappa$ is smooth and surjective. It suffices to show that 
$$
\kappa^*((\tau_1^*F)\ast (\tau_2^*F))
$$ 
is a shifted perverse sheaf, the intermediate extension from $\kappa^{-1}(\tilde U\times_{\cL(M)}\tilde U)_0$.

Let $\mu: \cW\to U\times_{\cL(M)} U$ be the map sending a point of $\cW$ to $(R,S,L)$. Applying (\ref{formula_explicit_F_appendix}) several times as in Step 3, we learn that there is a local system of rank one and order two, say $\cI$ on $\cW$ such that
$$
\kappa^*((\tau_1^*F)\ast (\tau_2^*F))\,\iso\, \cI\otimes \mu^*\pi_Y^*F
$$
Since $F$ is an irreducible perverse sheaf, our assertion follows.
\end{Prf}

\medskip

 Thus, Theorem~\ref{Th_appendix_CIO} is proved. 

\medskip\noindent
3.5 Now given $k$-points $N^0,L^0\in\wt\cL(M)$, let $F_{N^0,L^0}\in \D(H)$ be the $*$-restriction of $F$ under $(N^0,L^0)\times\id: 
H\hook{} \tilde Y\times H$. Define the functor $\cF_{N^0,L^0}: \D\!\cH_L\to\D\!\cH_N$ by 
$$
\cF_{N^0,L^0}(K)=F_{N^0,L^0}\ast K
$$
To see that it preserves perversity we can pick $S^0\in \wt\cL(M)$ with $N\cap S=L\cap S=0$ and use $\cF_{N^0,L^0}=\cF_{N^0,S^0}\comp \cF_{S^0,L^0}$. This reduces the question to the case $N\cap L=0$, in the latter case $\cF_{N^0,L^0}$ is nothing but the Fourier transform. 

 By Theorem~\ref{Th_appendix_CIO}, for $N^0,L^0,R^0\in\wt\cL(M)$ the diagram is canonically 2-commutative
$$
\begin{array}{ccc}
\D\!\cH_L & \;\toup{\cF_{R^0,L^0}} & \;\;\;\;\D\!\cH_R\\
 & \searrow\lefteqn{\scriptstyle \cF_{N^0, L^0}} & \;\;\;\;\downarrow\lefteqn{\scriptstyle \cF_{N^0, R^0}}\\
&& \;\;\;\;\D\!\cH_N
\end{array}
$$

\bigskip\noindent
3.6 {\scshape Nonramified Weil category}

\medskip\noindent
For a $k$-point $L^0\in \wt\cL(M)$ let $i_{L^0}: 
\wt\cL(M)\to \wt\cL(M)\times\wt\cL(M)\times H$ be the map sending $N^0$ to $(N^0, L^0, 0)$. We get a functor $\cF_{L^0}: \D\!\cH_L\to \D(\wt\cL(M))$ sending $K$ to the complex 
$$
i_{L^0}^*(F\ast \pr_3^*K)\otimes(\Qlb[1](\frac{1}{2}))^{\dim\cL(M)-2d-1}
$$

For any $k$-points $L^0,N^0\in\wt\cL(M)$ the diagram commutes
\begin{equation}
\label{diag_cF_transitivity}
\begin{array}{ccc}
\D\!\cH_L & \toup{\cF_{L^0}} & \;\;\D(\wt\cL(M))\\
& \searrow\lefteqn{\scriptstyle \cF_{L^0,N^0}} & \;\;\uparrow\lefteqn{\scriptstyle \cF_{N^0}}\\
  && \;\;\D\!\cH_N
\end{array}
\end{equation}
One checks that $\cF_{L^0}$ is exact for the perverse t-structure. 
\begin{Def} \select{The non-ramified Weil category} $W(\wt\cL(M))$ is the essential image of $\cF_{L_0}: \cH_L\to \P(\wt\cL(M))$. 
This is a full subcategory in $\P(\wt\cL(M))$ independent of $L^0$, because (\ref{diag_cF_transitivity}) is commutative.
\end{Def}

 The group $G$ acts naturally on $\wt\cL(M)$, hence also on $\P(\wt\cL(M))$. This action preserves the full subcategory $W(\wt\cL(M))$. 
 
 At the classical level, for $L\in\cL(M)$ the $G$-respresentation $\cH_L\,\iso\, \cH_{L, odd}\oplus \cH_{L, even}$ is a direct sum of two irreducible ones consisting of odd and even functions respectively. 
The category $W(\wt\cL(M))$ is a geometric analog of the space $\cH_{L, even}$. The geometric analog of the whole Weil representation $\cH_L$ is as follows.

\begin{Def} Let $\act_l: \bar C\times H\to \wt\cL(M)\times H$ be the map sending $(L^0,h, \bar l\in\bar L)$ to $(L^0, \bar l h)$. A perverse sheaf $K\in\P(\wt\cL(M)\times H)$ is $(\bar C, \chi_{\bar C})$-\select{equivariant} if it is equipped with an isomorphism 
$$
\act_l^*K\,\iso\, \pr^*K\otimes \pr_1^*\chi_{\bar C}
$$ 
satisfying the usual associativity property, and whose restiction to the unit section is the identity.

\select{The complete Weil category} $W(M)$ is the category of pairs $(K,\sigma)$, where $K\in \P(\wt\cL(M)\times H)$ is a $(\bar C, \chi_{\bar C})$-equivariant perverse sheaf, and 
$$
\sigma: F\ast\pr_{23}^*K\,\iso\, \pr_{13}^*K
$$
is an isomorphism for the projections $\pr_{13},\pr_{23}: \wt\cL(M)\times\wt\cL(M)\times H\to \wt\cL(M)\times H$. The map $\sigma$ must be compatible with the associativity constraint and the unit section constraint of $F$. 
\end{Def}

 The group $G$ acts on $\wt\cL(M)\times H$ sending $(g\in G, L^0,h)$ to $(gL^0, gh)$. This action extends to an action of $G$ on the category $W(M)$. 
 
\bigskip\medskip

\centerline{\scshape 4. Compatibility property}

\medskip\noindent
4.1 In this section we establish the following additional property of the canonical interwining operators. 
Let $V\subset M$ be an isotropic subspace, $V^{\perp}\subset M$ its orthogonal complement. Let 
$\cL(M)_V\subset \cL(M)$ be the open subscheme of $L\in\cL(M)$ such that $L\cap V=0$. Set $M_0=V^{\perp}/V$. We have a map 
$p_V: \cL(M)_V\to \cL(M_0)$ sending $L$ to $L_V:=L\cap V^{\perp}$. 

 Write $Y=\cL(M)\times\cL(M)$ and $Y_V=\cL(M)_V\times\cL(M)_V$. The gerb $\wt Y$ is defined as in Section~3.2, write $\wt Y_V$ for its restriction to $Y_V$. Set  $Y_0=\cL(M_0)\times\cL(M_0)$, we have the corresponding gerb $\wt Y_0$ defined as in Section~3.2. We extend the map $p_V\times p_V$ to a map
$$
\pi_V: \wt Y_V\to \wt Y_0
$$
sending $(L_1,L_2, \cB,\,  \cB^2\,\iso\, \det L_1\otimes\det L_2)$ to 
$$
(L_{1,V},L_{2,V}, \cB_0, \;\cB_0^2\,\iso\, \det L_{1,V}\otimes\det L_{2,V})
$$ 
Here $L_{i,V}=L_i\cap V^{\perp}$ and $\cB_0=\cB\otimes\det V$. 
 We used the exact sequences 
$$
0\to L_{i,V}\to L_i\to M/V^{\perp}\to 0
$$ 
yielding canonical ($\ZZ/2\ZZ$-graded) isomorphisms $\det L_{i,V}\otimes\det V^*\,\iso\, \det L_i$. 

Write $H_0=M_0\oplus k$ for the Heisenberg group of $M_0$. For $L\in\cL(M)_V$ we have the categories $\cH_L$ and $\cH_{L_V}$ of certain perverse sheaves on $H$ and $H_0$ respectively. To such $L$ we associate a transition functor $T^L: {\cH_{L_V}}\to \cH_L$ which will be fully faithful and exact for the perverse t-structures. 

 Write for brevity $H^V=V^{\perp}\times \A^1$. First, at the level of functions, given $f\in \cH_{L_V}$ consider it as a function on $H^V$ via the composition $H^V\toup{\alpha_V} H_0\toup{f}\Qlb$, where $\alpha_V$ sends $(v,a)$ to $(v\mod V, a)$. Then there is a unique $f_1\in \cH_L$ such that $f_1(m)=q^{\dim V}f(m)$ for all $m\in H^V$. We use the property $V^{\perp}+L=M$. We set 
\begin{equation}
\label{def_T^L_classical}
(T^L)(f)=f_1
\end{equation}
The image of $T^L$ is 
$$
\{f_1\in \cH_L\mid f(h(v,0))=f(h),\;\; h\in H, v\in V\}
$$ 
Note that $H^V\subset H$ is a subgroup, and $V=\{(v,0)\in H^V\mid v\in V\}\subset H^V$ is a normal subgroup lying in the center of $H^V$. The operator $T^L:\cH_{L_V}\to \cH_L$ commutes with the action of $H^V$. It is understood that on $\cH_{L_V}$ this group acts via its quotient $H^V\toup{\alpha_V} H_0$.  

 On the geometric level, consider the map $s: L\times H^V\to H$ sending $(l, (v,a))$ to the product in the Heisenberg group $(l,0)(v,a)\in H$. Note that $s$ is smooth and surjective, an affine fibration of rank $\dim L_V$. Given $K\in \cH_{L_V}$ there is a (defined up to a unique isomorphism) perverse sheaf $T^LK\in \cH_L$ equipped with
$$
s^*(T^LK)\otimes(\Qlb[1](\frac{1}{2}))^{\dim L_V}\,\iso\,
\Qlb\boxtimes\alpha_V^*K\otimes(\Qlb[1](\frac{1}{2}))^{\dim V+\dim L}
$$

 The \select{compatibility property} of the canonical interwining operators is as follows. 

\begin{Pp} 
\label{Pp_appendix_compatibility_property_CIO}
Let $(L,N,\cB)\in \wt Y_V$, write $(L_V,N_V,\cB_0)$ for the image of $(L,N,\cB)$ under $\pi_V$. Write $\cF_{N^0,L^0}: \cH_L\to\cH_N$ and $\cF_{N^0_V, L^0_V}: \cH_{L_V}\to \cH_{N_V}$ for the corresponding functors defined as in Section~3.5. Then the diagram of categories is canonically 2-commutative
$$
\begin{array}{ccccccc}
\cH_{L_V} &&& \toup{T^L} &&& \cH_L\\
\downarrow\lefteqn{\scriptstyle \cF_{N_V^0, L_V^0}} &&&&&& \downarrow\lefteqn{\scriptstyle \cF_{N^0, L^0}}\\
\cH_{N_V} &&& \toup{T^N} &&& \cH_N
\end{array}
$$
One may also replace $\cH$ by $\D\!\cH$ in the above diagram. 
\end{Pp}

\medskip\noindent
4.2 First, we realize the functors $T^L$ by a universal kernel, namely, we define a perverse sheaf $T$ on $\cL(M)_V\times H\times H_0$ as follows. 
 
 Remind the vector bundle $\bar C\to \cL(M)$, its fibre over $L$ is $\bar L=L\times\A^1$. Write $\bar C_V$ for 
the restriction of $\bar C$ to the open subscheme $\cL(M)_V$. 
We have a closed immersion 
$$
i_0: \bar C_V\times H^V\to \cL(M)_V\times H\times H_0
$$ 
sending $(\bar l\in\bar L, u\in H^V)$ to $(L, \bar l u, \alpha_V(u))$, where the product $\bar l u$ is taken in $H$. The perverse sheaf $T$ is defined by
$$
T=(i_0)_!\pr_1^*\chi_{\tilde C}\otimes(\Qlb[1](\frac{1}{2}))^{\dim\bar C+\dim V+\dim H_0},
$$
here $\pr_1: \bar C_V\times H^V\to \bar C_V$ is the projection, and $\chi_{\bar C}$ was defined in 3.1. 

 For $L\in \cL(M)_V$ let $T_L$ be the $*$-restriction of $T$ under $(L,\id): H\times H_0\to \cL(M)_V\times H\times H_0$. Define $T^L: \D\!\cH_{L_V}\to \D\!\cH_L$ by
\begin{equation}
\label{def_T^L_geom}
T^L(K)\,\iso\, \pr_{1!}(T_L\otimes \pr_2^*K)\otimes(\Qlb[1](\frac{1}{2}))^{\dim V-d-\dim\cL(M)}
\end{equation}
for the diagram of projections $H\getsup{\pr_1} H\times H_0\toup{\pr_2} H_0$. It is exact for the perverse t-structures. 
 
 The sheaf $T$ has the following properties. At the level of functions, the corresponding function $T_L: H\times H_0\to\Qlb$ satisfies 
$$
T_L(\bar l h, \bar l_0 h_0)=\chi_L(\bar l)\chi_{L_V}(\bar l_0)^{-1} T_L(h,h_0),\;\;\; \bar l\in\bar L, \bar l_0\in\bar L_V
$$
The geometric analog is as follows. Let $^0\bar C\to \cL(M)_V$ be the vector bundle, whose fibre over $L\in\cL(M)_V$ is $\bar L\times \bar L_V$. Consider the diagram 
$$
\cL(M)_V\times H\times H_0\getsup{\pr^V}  {^0\bar C}\times H\times H_0\toup{\act_{lr}^V} \cL(M)_V\times H\times H_0,
$$
where $\pr^V$ is the projection, and $\act_{lr}^V$
sends 
$$
(L\in\cL(M)_V, \bar l\in \bar L, \bar l_0\in\bar L_V, h\in H,h_0\in H_0)
$$ 
to $(L, \bar l h, \bar l_0 h_0)$. Let $^0p: {^0\bar C}\to\A^1$ be the map sending 
$$
(L\in\cL(M)_V, \bar l\in \bar L, \bar l_0\in\bar L_V)
$$ 
to $p(\bar l)-p(\bar l_0)$. Here $p: \bar L\to\A^1$ and $p:\bar L_V\to\A^1$ are the projections on the center. Set $^0\chi=(^0p)^*\cL_{\psi}$. Then $T$ is \select{$\act_{lr}^V$-equivariant}, that is, it admits an isomorphism
$$
(\act_{lr}^V)^*T\,\iso\, (\pr^V)^*T\otimes \pr_1^*(^0\chi),
$$
satisfying the usual associativity property, and its restriction to the unit section is the identity. 

\medskip\noindent
4.3 We will prove a geometric version of the equality 
(up to an explicit power of $q$)
$$
\int_{u\in H} F_{N^0,L^0}(hu^{-1})T_L(u, h_0)du=
\int_{v\in H_0} T_N(h, v)F_{N_V^0, L^0_V}(vh_0^{-1})dv
$$
for $h\in H, h_0\in H_0$. Here $(N^0,L^0)\in \tilde Y_V$ and 
$$
(N^0_V, L^0_V)=\pi_V(N^0,L^0)
$$ 

  Write $\inv: H\iso H$ for the map sending $h$ to $h^{-1}$, set $^{\inv}F=(\id\times\inv)^*F$ for $\id\times\inv: \wt Y\times H\to \wt Y\times H$. For $i=1,2$ write $p_i: \wt Y_V\to\cL(M)_V$ for the projection on the $i$-th factor.  Let $q_0$ denote the composition
$$
\tilde Y_V\times H\times H_0\toup{\pr_{13}} \tilde Y_V\times H_0\toup{\pi_V\times\id} \tilde Y_0\times H_0
$$   
Proposition~\ref{Pp_appendix_compatibility_property_CIO} is an immediate consequence of the following.
\begin{Lm}  
\label{Lm_appendix_compatibility_property}
There is a canonical isomorphism over $\wt Y_V\times H\times H_0$
$$
(\pr_{12}^*F)\ast_H (p_2\times\id)^*T\,\iso\, (q_0^*(^{\inv}F))\ast_{H_0} (p_1\times\id)^*T
$$
where $\pr_{12}: \tilde Y_V\times H\times H_0\to \tilde Y_V\times H$ and $\, p_1\times\id, p_2\times\id: \tilde Y_V\times H\times H_0\to \cL(M)_V\times H\times H_0$. 
\end{Lm}

Let $i_V: H^V\hook{} H$ be the natural closed immersion.
It is elementary to check that Lemma~\ref{Lm_appendix_compatibility_property} is equivalent to the following.  
\begin{Lm}  
\label{Lm_long}
There is a canonical isomorphism of (shifted) perverse sheaves 
\begin{equation}
\label{iso_compatibility_F'}
(\id\times\alpha_V)_!i_V^*F \,\iso\, (\pi_V\times\id)^*F\otimes(\Qlb[1](\frac{1}{2}))^{\dimrel(\pi_V)+\dim V}
\end{equation}
for the diagram
$$
\begin{array}{ccc}
& \wt Y_V\times H^V & \toup{i_V} \wt Y_V\times H\\
& \downarrow\lefteqn{\scriptstyle \id\times\alpha_V}\\
\wt Y_0\times H_0 \getsup{\pi_V\times\id} & \wt Y_V\times H_0
\end{array}
$$
\end{Lm}  
\begin{Prf} 
Write $U(M_0)$ for the scheme $U$ constructed out of the symplectic space $M_0$, it classifies pairs of lagrangian subspaces in $M_0$ that do not intersect. 
We have a 2-commutative diagram
$$
\begin{array}{ccccc}
U(M_0)\times_{\cL(M_0)} U(M_0) & \getsup{\pi_W} & W_V & \hook{i_W} & U\times_{\cL(M)} U\\
\downarrow\lefteqn{\scriptstyle \pi_{Y_0}} && \downarrow\lefteqn{\scriptstyle \pi_{Y,V}} & \;\;\;\swarrow\lefteqn{\scriptstyle \pi_Y}\\
\wt Y_0 & \getsup{\pi_V} & \wt Y_V
\end{array}
$$
where the square is cartesian thus defining $W_V, \pi_W$, and $\pi_{Y,V}$. The map $i_W$ is a locally closed immersion. Write a point of $W_V$ as a triple $(N,R,L)\in \cL(M)$ such that $N,L\in \cL(M)_V$, $V\subset R\subset V^{\perp}$, and $N\cap R=R\cap L=0$. The map $\pi_W$ sends $(N,R,L)$ to $(N_V, R_V, L_V)$ with $R_V=R/V$. 

 Let us establish the isomorphism (\ref{iso_compatibility_F'}) after restriction under $\pi_{Y,V}\times\alpha_V: W_V\times H^V\to \tilde Y_V\times H_0$. We first give the argument at the level of functions and then check that it holds through in the geometric setting.
 
 Consider a point of $W_V$ given by a triple $(N,R,L)\in\cL(M)$, so $N,L\in \cL(M)_V$, $V\subset R\subset V^{\perp}$, and $N\cap R=R\cap L=0$. We have $V^{\perp}=R\oplus L_V$. Let $h\in H^V$, write $h=(r,a)(l_1,0)$ for uniquely defined $r\in R, l_1\in L_V, a\in k$.  Write $(N^0,L^0)\in\wt Y_V$ for the image of $(N,R,L)$ under $\pi_{Y,V}$. Using (\ref{formula_explicit_F_appendix}), we get
\begin{multline*}
\int_{v\in V} F_{N^0, L^0}(h(v,0))dv=q^{\dim\cL(M)-\frac{d+1}{2}}\int_{v\in V, u\in H} \tilde F_{N,R}(u)\tilde F_{R,L}(u^{-1}h(v,0))dvdu=\\
q^{\dim\cL(M)+\frac{d+1}{2}}\int_{v\in V, u\in H/\tilde R}  \tilde F_{N,R}(u)\tilde F_{R,L}(u^{-1}(r,a)(v,0))dvdu=\\
q^{\dim\cL(M)+\frac{d+1}{2}}\int_{v\in V, l\in L}  \tilde F_{N,R}(l,0)\tilde F_{R,L}((-l,0)(r,a)(v,0))dvdl
\end{multline*}
Since $(-l,0)(r+v,a)=(r+v, a+\omega\<r+v,l\>)(-l,0)$, the latter expression equals 
$$
q^{-\frac{d}{2}}\int_{v\in V, l\in L}  \tilde F_{N,R}(l,0)\psi(a+\omega\<r+v,l\>)dvdl=
q^{\dim V-\frac{d}{2}}\int_{l\in L_V} 
\tilde F_{N,R}(l,0)\psi(a+\omega\<r,l\>)dl
$$
For $l\in L_V$ we get $\tilde F_{N,R}(l,0)=q^{\dim\cL(M_0)-\dim\cL(M)-\dim V}\tilde F_{N_V,R_V}(l,0)$. Indeed, since $V^{\perp}=R\oplus N_V$, there are unique $r_1\in R, n_1\in N_V$ such that $l=n_1+r_1$. For $\bar r_1=r_1\!\mod V\in M_0$ we get 
\begin{multline*}
\tilde F_{N,R}(l,0)=q^{-\dim\cL(M)-\frac{2d+1}{2}}
\chi_{NR}(l,0)=q^{-\dim\cL(M)-\frac{2d+1}{2}}\psi(\frac{1}{2}\omega\<r_1,n_1\>)=\\ q^{-\dim\cL(M)-\frac{2d+1}{2}}\chi_{N_VR_V}(\bar r_1+n_1,0)=q^{\dim\cL(M_0)-\dim\cL(M)-\dim V}\tilde F_{N_V,R_V}(l,0)
\end{multline*}
Further, we claim that 
$$
\tilde F_{R_V,L_V}((-l,0)\alpha_V(h))=q^{-\dim\cL(M_0)-\frac{\dim H_0}{2}}\psi(a+\omega\<r,l\>)
$$ 
This follows from definition (\ref{def_tilde_F_U}) of $\tilde F_U$ and the formula $(-l,0)(r,a)=(r, a+\omega\<r,l\>)(-l,0)$. 
 
  Combinig the above we get
\begin{multline*}
\int_{v\in V} F_{N^0, L^0}(h(v,0))dv=q^c
\int_{l\in L_V} \tilde F_{N_V,R_V}(l,0)\tilde F_{R_V,L_V}((-l,0)\alpha_V(h))dl=\\
q^{c+\dim V-d-1}
\int_{u\in H_0} \tilde F_{N_V,R_V}(u)\tilde F_{R_V,L_V}(u^{-1}\alpha_V(h))du
\end{multline*}
with $c=\frac{\dim H_0-d}{2}+2\dim\cL(M_0)-\dim\cL(M)$. 
By (\ref{formula_explicit_F_appendix}), the latter expression identifies with $F_{N_V^0, L_V^0}(h)$ up to an explicit power of $q$. 

 The argument holds through in the geometric setting yielding the desired isomorphism $\gamma$ over $W_V\times H^V$. For any point $(N_V, L_V\cB_0)\in \tilde Y_0$ such that $N_V\ne L_V$ the fibre of $\pi_{Y_0}$ over this point is geometrically connected. So, for $\dim V<d$ the isomorphism $\gamma$ descends  to a uniquely defined isomorphism (\ref{iso_compatibility_F'}). The case $\dim V=d$ is easier and is left to the reader.
\end{Prf}

\medskip

\begin{Rem} Let $i_H:\Spec k\hook{} H$ denote the zero section. Arguing as in Lemma~\ref{Lm_long}, for the map $\id\times i_H: \wt Y\to \wt Y\times H$ one gets a canonical isomorphism 
$$
(\id\times i_H)^*F\,\iso\, (S_{M,g}\otimes\cE^{d}\oplus S_{M,s}\otimes\cE^{d-1})\otimes(\Qlb[1](\frac{1}{2}))^{\dim H},
$$
it will not be used in this paper.
\end{Rem}

\medskip\noindent
4.4 The functors $T^L$ satisfy the following transitivity property. Assume that $V_1\subset V$ is another isotropic subspace in $M$. Let $M_1=V_1^{\perp}/V_1$ 
and $H_1=M_1\times\A^1$ be the corresponding Heisenberg group. Then for $L\in \cL(M)_V$ we also have $L_{V_1}:=L\cap V_1^{\perp}$ and the category $\cH_{L_{V_1}}$ of certain perverse sheaves on $H_1$. Then the diagram is canonically 2-commutative
$$
\begin{array}{ccc}
\cH_{L_V} & \toup{T^{L_{V_1}}} & \cH_{L_{V_1}}\\
& \searrow\lefteqn{\scriptstyle T^L} & \downarrow\lefteqn{\scriptstyle T^L}\\
&& \cH_L
\end{array}
$$

\smallskip\noindent
4.5 We will need also one more compatibility property of the canonical interwining operators. Let $V\subset V^{\perp}\subset M$ be as in 4.1. Write $i_{0,V}: \cL(M_0)\to \cL(M)$ for the closed immersion sending $L_0$ to the preimage of $L_0$ under $V^{\perp}\to V^{\perp}/V$.

 For $L\in\cL(M)$ with $V\subset L$ set $L_V=L/V\in \cL(M_0)$. Let $(\cL(M_0)\times\cL(M)_V)^{\tilde{}}$ denote the restriction of the gerb $\tilde Y$ under 
$$
\cL(M_0)\times\cL(M)_V\toup{i_{0,V}\times\id} \cL(M)\times \cL(M)_V\subset Y
$$ 
Define $\pi_{0,V}: (\cL(M_0)\times\cL(M)_V)^{\tilde{}}\to \tilde Y_0$ as the map sending $(L,N,\cB, \cB^2\,\iso\, \det L\otimes\det N)$ to 
$$
(L_V, N_V, \cB, \cB^2\,\iso\, \det L_V\otimes\det N_V)
$$ 
Here $L\in\cL(M)$ with $V\subset L$. We have used the canonical $\ZZ/2\ZZ$-graded isomorphism 
$\det L\otimes\det N\,\iso\, \det L_V\otimes\det N_V$. 

 Remind the closed immersion $i_V: H^V\hook{} H$. 
For $L\in\cL(M)$ with $V\subset L$ define the transition functor $T^L: \cH_{L_V}\to \cH_L$ by
$$
T^L(K)=i_{V !}\alpha_V^*K\otimes(\Qlb[1](\frac{1}{2}))^{\dim V}
$$

 The proof of the following is similar to that of Proposition~\ref{Pp_appendix_compatibility_property_CIO} and is left to the reader. 
\begin{Pp}  
\label{Pp_second_compatibility}
Let $(L,N,\cB)\in (\cL(M_0)\times\cL(M)_V)^{\tilde{}}$, let $(L_V, N_V, \cB)$ denote its image under $\pi_{0,V}$. Write $\cF_{N^0, L^0}: \cH_L\to \cH_N$ and $\cF_{N_V^0, L_V^0}$ for the corresponding functors defined as in Section~3.5. Then the diagram of categories is canonically 2-commutative
$$
\begin{array}{ccc}
\cH_{L_V} \;\;\;\;\;\;& \toup{T^L} & \;\;\cH_L\\
\downarrow\lefteqn{\scriptstyle \cF_{N^0_V, L^0_V}} \;\;\;\;\;\;&& \;\;\downarrow\lefteqn{\scriptstyle \cF_{N^0, L^0}}\\
\cH_{N_V} \;\;\;\;\;\;&\toup{T^N} & \;\;\cH_N
\end{array}
$$
One may also replace $\cH$ by $\D\!\cH$ in the above diagram. \QED
\end{Pp}

\bigskip\medskip

\centerline{\scshape 5. Discrete lagrangian lattices and the metaplectic group}

\bigskip\noindent
5.1 Set $\cO=k[[t]]\subset F=k((t))$. Denote by $\Omega$ the completed module of relative differentials of $\cO$ over $k$. Let $M$ be a free $\cO$-module of rank $2d$ with symplectic form $\wedge^2 M\to\Omega$. Write $G$ for the group scheme over $\Spec\cO$ of automorphisms of $M$ preserving the symplectic form. Consider the Tate space $M(F)$ (cf. \cite{BD}, 4.2.13 for the definition), it is equipped with the symplectic form $(m_1,m_2)\mapsto \Res\omega\<m_1,m_2\>$. 

 For a $k$-subspace $L\subset M(F)$ write 
$$
L^{\perp}=\{m\in M(F)\mid \Res\omega\<m,l\>=0\;\;\mbox{for all}\;\; l\in L\}
$$ 
For two $k$-subspaces $L_1,L_2\subset M$ we get $(L_1+L_2)^{\perp}=L_1^{\perp}\cap L_2^{\perp}$. 
For a finite-dimensional symplectic $k$-vector space $U$ write $\cL(U)$ for the variety of lagrangian subspaces in $U$. 

 As in \select{loc.cit}, we say that an $\cO$-submodule $R\subset M(F)$ is a \select{c-lattice} if $M(-N)\subset R\subset M(N)$ for some integer $N$. A \select{lagrangian d-lattice} in $M(F)$ is a $k$-vector subspace $L\subset M(F)$ such that $L^{\perp}=L$ and there exists a $c$-lattice $R$ with $R\cap L=0$. Note that the condition $R\cap L=0$ implies $R^{\perp}+L=M(F)$. 
Let $\cL_d(M(F))$ denote the set of lagrangian d-lattices in $M(F)$. 
  
 For a given c-lattice $R\subset M(F)$ write 
$$
\cL_d(M(F))_R=\{L\in\cL_d(M(F))\mid L\cap R=0\}
$$ 
If $R$ is a c-lattice in $M(F)$ with $R\subset R^{\perp}$ then $\cL_d(M(F))_R$ is a naturally a $k$-scheme (not of finite type over $k$). Indeed, for each c-lattice $R_1\subset R$ we have the variety 
$$
\cL(R_1^{\perp}/R_1)_R:=\{L_1\in \cL(R_1^{\perp}/R_1)\mid L_1\cap R/R_1=0\}
$$
For $R_2\subset R_1\subset R$ we get a map $p_{R_2,R_1}:  \cL(R_2^{\perp}/R_2)_R \to \cL(R_1^{\perp}/R_1)_R$ sending $L_2$ to 
$$
L_1=(L_2\cap (R_1^{\perp}/R_2))+R_1
$$ 
The map $p_{R_2,R_1}$ is a composition of two affine fibrations of constant rank. Then $\cL_d(M(F))_R$ is the inverse limit of $\cL(R_1^{\perp}/R_1)_R$ over the partially ordered set of c-lattices $R_1\subset R$. 

  If $R'\subset R$ is another c-lattice then $\cL_d(M(F))_R\subset \cL_d(M(F))_{R'}$ is an open immersion (as it is an open immersion on each term of the projective system). So, $\cL_d(M(F))$ is a $k$-scheme that can be seen as the inductive limit of $\cL_d(M(F))_R$. 

 Let us define the categories $\P(\cL_d(M(F)))$ and 
 $\P_{G(\cO)}(\cL_d(M(F)))$ of perverse sheaves and $G(\cO)$-equivariant perverse sheaves on $\cL_d(M(F))$. 
 
 For $r\ge 0$ set 
$$
_r\cL_d(M(F))=\cL_d(M(F))_{M(-r)},
$$ 
the group $G(\cO)$ acts on $_r\cL_d(M(F))$ naturally. First, define the category $\D_{G(\cO)}(_r\cL_d(M(F)))$ as follows. 

 For $N+r\ge 0$ set $_{N,r}M=t^{-N}M/t^rM$. 
For $N\ge r\ge 0$ the action of $G(\cO)$ on $_r\cL(_{N,N}M):=\cL(_{N,N}M)_{M(-r)}$ factors through $G(\cO/t^{2N})$. For $r_1\ge 2N$ the kernel 
$$
\Ker(G(\cO/t^{r_1}))\to G(\cO/t^{2N}))
$$ 
is unipotent, so that we have an equivalence (exact for the perverse t-structures) 
$$
\D_{G(\cO/t^{2N})}(_r\cL(_{N,N}M))\,\iso\, \D_{G(\cO/t^{r_1})}(_r\cL(_{N,N}M))
$$
Define $\D_{G(\cO)}(_r\cL(_{N,N}M))$ as $\D_{G(\cO/t^{r_1})}(_r\cL(_{N,N}M))$ for any $r_1\ge 2N$. It is equipped with the perverse t-structure. 

For $N_1\ge N\ge r\ge 0$ the fibres of the above projection 
$$
p: {_r\cL(_{N_1,N_1}M)}\to {_r\cL(_{N,N}M)}
$$ 
are isomorphic to affine spaces of fixed dimension, and $p$ is smooth and surjective. Hence, this map yields  transition functors (exact for the perverse t-structures and fully faithful embeddings)
$$
\D_{G(\cO)}(_r\cL(_{N,N}M))\to 
\D_{G(\cO)}(_r\cL(_{N_1,N_1}M))
$$
and 
$$
\D(_r\cL(_{N,N}M))\to 
\D(_r\cL(_{N_1,N_1}M))
$$
We define $\D_{G(\cO)}(_r\cL_d(M(F)))$ as the inductive 2-limit of $\D_{G(\cO)}(_r\cL(_{N,N}M))$ as $N$ goes to plus infinity. The category $\D(_r\cL_d(M(F)))$ is defined similarly. Both they are equipped with perverse t-structures. 

 If $N_1\ge N\ge r_1\ge r\ge 0$ we have a diagram
$$
\begin{array}{ccc}
_r\cL(_{N_1,N_1}M) & \toup{p} & {_r\cL(_{N,N}M)}\\
\downarrow\lefteqn{\scriptstyle j} && \downarrow\lefteqn{\scriptstyle j}\\
_{r_1}\cL(_{N_1,N_1}M) & \toup{p} & {_{r_1}\cL(_{N,N}M)},
\end{array}
$$
where $j$ are natural open immersions. The restriction functors $j^*: \D_{G(\cO)}(_{r_1}\cL(_{N,N}M))\to \D_{G(\cO)}(_r\cL(_{N,N}M))$ yield (in the limit as $N$ goes to plus infinity) the functors 
$$
j_{r_1,r}^*: \D_{G(\cO)}(_{r_1}\cL_d(M(F)))\to \D_{G(\cO)}(_r\cL_d(M(F)))
$$
of restriction with respect to the open immersion $j_{r_1,r}: {_r\cL_d(M(F))}\hook{} {_{r_1}\cL_d(M(F))}$. 
Define
$\D_{G(\cO)}(\cL_d(M(F)))$ as the projective 2-limit of 
$$
\D_{G(\cO)}(_r\cL_d(M(F)))
$$ as $r$ goes to plus infinity. Similarly, $\P_{G(\cO)}(\cL_d(M(F)))$ is defined as the projective 2-limit of $\P_{G(\cO)}(_r\cL_d(M(F)))$. Along the same lines, one defines the categories $\P(\cL_d(M(F)))$ and $\D(\cL_d(M(F)))$.
 
\medskip\noindent
5.2 {\scshape Relative determinant\  } 
For a pair of c-lattices $M_1,M_2$ in $M(F)$ define the relative determinant $\det(M_1: M_2)$ as the following
$\ZZ/2\ZZ$-graded 1-dimensional $k$-vector space.
If $R$ is a c-lattice in $M(F)$ such that $R\subset M_1\cap M_2$ then 
$$
\det(M_1:M_2)\,\iso\, \det(M_1/R)\otimes\det(M_2/R)^{-1},
$$ 
it is defined up to a unique isomorphism. 
 
 Write $\Gr_G$ for the affine grassmanian $G(F)/G(\cO)$ of $G$ (cf.~\cite{BD}, Section~4.5). For $R\in\Gr_G, L\in \cL_d(M(F))$ define the relative determinant $\det(R:L)$ as the following ($\ZZ/2\ZZ$-graded purely of degree zero) 1-dimensional vector space. Pick a c-lattice $R_1\subset R$ such that $R_1\cap L=0$. Then in $R_1^{\perp}/R_1$ one gets two lagrangian subspaces $R/R_1$ and $L_{R_1}:=L\cap R_1^{\perp}$. Set 
$$
\det(R:L)=\det(R/R_1)\otimes\det(L_{R_1})
$$ 
If $R_2\subset R_1$ is another c-lattice then the exact sequence 
$$
0\to L_{R_1}\to L\cap R_2^{\perp}\to R_2^{\perp}/R_1^{\perp}\to 0
$$
yields a canonical $\ZZ/2\ZZ$-graded isomorphism
\begin{multline*}
\det(R/R_2)\otimes\det(L_{R_2})\,\iso\,
\det(R_1/R_2)\otimes\det (R/R_1)\otimes \det(L_{R_1})\otimes \det(R_2^{\perp}/R_1^{\perp})\,\iso\\ 
\det(R/R_1)\otimes\det(L_{R_1})
\end{multline*} 
So, $\det(R:L)$ is a $\ZZ/2\ZZ$-graded line defined up to a unique isomorphism. Another way to say is as follows. Consider the complex $R\oplus L\toup{s} M(F)$ placed in cohomological degrees 0 and 1, where $s(r,l)=r+l$. It has finite-dimensional cohomologies and 
$$
\det(R:L)=\det(R\oplus L\toup{s} M(F))
$$

 For $g\in G(F)$ we have canonically 
$$
\det(gR: gL)\,\iso\, \det(R:L)
$$ 
For $R_1,R_2\in\Gr_G$, $L\in\cL_d(M(F))$ we have canonically
$$
\det(R_1: L)\,\iso\, \det(R_1: R_2)\otimes\det(R_2: L)
$$
 
\noindent 
5.3 Write $\cA_d$ for the line bundle on $\cL_d(M(F))$ with fibre $\det(M: L)$ at $L\in \cL_d(M(F))$. Clearly, $\cA_d$ is $G(\cO)$-equivariant, so we may see $\cA_d$ as the line bundle on the stack quotient $\cL_d(M(F))/G(\cO)$. Let $\tilde\cL_d(M(F))$ denote the $\mu_2$-gerb of square roots of $\cA_d$. 

 The categories of the corresponding perverse sheaves $\P_{G(\cO)}(\tilde\cL_d(M(F))$ and $\P(\tilde\cL_d(M(F))$ are defined as above. Namely, first for $r\ge 0$ define
$$
\D_{G(\cO)}(_r\tilde\cL_d(M(F)))
$$ 
as follows. For $N\ge r$ take $r_1\ge 2N$ and consider the stack quotient $_r\cL(_{N,N}M)/G(\cO/t^{r_1})$. We have the line bundle, say $\cA_N$ on this stack whose fibre at $L$ is $\det(M/M(-N))\otimes \det L$. Here $L\subset {_{N,N}M}$ is a Lagrangian subspace such that $L\cap (M(-r)/M(-N))=0$. Write 
$$
(_r\cL(_{N,N}M)/G(\cO/t^{r_1}))\tilde{}
$$
for the gerb of square roots of this line bundle. Let 
$\D_{G(\cO)}(_r\tilde\cL(_{N,N}M))$ denote the category 
$$
\D((_r\cL(_{N,N}M)/G(\cO/t^{r_1}))\tilde{})
$$ 
for any $r_1\ge 2N$  (we have canonical equivalences exact for the perverse t-strucures between such categories for various $r_1$).

 Assume $N_1\ge N\ge r$ and $r_1\ge 2N_1$. For the projection
$$
p: {_r\cL(_{N_1,N_1}M)/G(\cO/t^{r_1})}\to  {_r\cL(_{N,N}M)/G(\cO/t^{r_1})}
$$
we have a canonical $\ZZ/2\ZZ$-graded isomorphism $p^*\cA_N\,\iso\, \cA_{N_1}$. This yields a transition map
$$
(_r\cL(_{N_1,N_1}M)/G(\cO/t^{r_1}))\tilde{}\to  (_r\cL(_{N,N}M)/G(\cO/t^{r_1}))\tilde{}
$$
The corresponding inverse image yields a transition functor 
\begin{equation}
\label{trans_functors_open_imm_N}
\D_{G(\cO)}(_r\tilde\cL(_{N,N}M))\to \D_{G(\cO)}(_r\tilde\cL(_{N_1,N_1}M))
\end{equation}
exact for the perverse t-structures (and a fully faithful embedding). We define $\D_{G(\cO)}(_r\tilde\cL_d(M(F)))$ as the inductive 2-limit of $\D_{G(\cO)}(_r\tilde\cL(_{N,N}M))$ as $N$ goes to plus infinity. 

 For $N\ge r'\ge r$ and $r_1\ge 2N$ we have an open immersion
$$
\tilde j: (_r\cL(_{N,N}M)/G(\cO/t^{r_1}))\tilde{}\subset 
(_{r'}\cL(_{N,N}M)/G(\cO/t^{r_1}))\tilde{}
$$
hence the $*$-restriction functors 
$$
\tilde j^*: \D_{G(\cO)}(_{r'}\tilde\cL(_{N,N}M))\to \D_{G(\cO)}(_r\tilde\cL(_{N,N}M))
$$  
compatible with the transition functors (\ref{trans_functors_open_imm_N}). Passing to the limit as $N$ goes to plus infinity, we get the functors 
$$
\tilde j^*_{r',r}: \D_{G(\cO)}(_{r'}\tilde\cL_d(M(F)))\to \D_{G(\cO)}(_r\tilde\cL_d(M(F)))
$$  
Define $\D_{G(\cO)}(\tilde\cL_d(M(F)))$ as the projective 2-limit of $\D_{G(\cO)}(_r\tilde\cL_d(M(F)))$ as $r$ goes to plus infinity, and similarly for $\P_{G(\cO)}(\tilde\cL_d(M(F)))$.

 Along the same lines one defines the categories $\P(\tilde\cL_d(M(F)))$ and $\D(\tilde\cL_d(M(F)))$.
 
\medskip\noindent
5.4 {\scshape Metaplectic group\ } Let $\cA_G$ be the line bundle on the ind-scheme $G(F)$ whose fibre at $g$ is $\det(M:gM)$. Write $\wt G(F)\to G(F)$ for the gerb of square roots of $\cA_G$. The stack $\wt G(F)$ has a structure of a group stack. The product map $m:\wt G(F)\times\wt G(F)\to\wt G(F)$ sends 
$$
(g_1,\cB_1, \sigma_1:\cB_1^2\,\iso\,\det(M: g_1M)), (g_2,\cB_2, \sigma_2: \cB_2^2\,\iso\, \det(M:g_2 M))
$$ 
to the collection $(g_1g_2,\cB, \sigma: \cB^2\,\iso\, \det(M:g_1g_2M)$, where $\cB=\cB_1\otimes\cB_2$ and 
$\sigma$ is the composition
\begin{multline*}
(\cB_1\otimes\cB_2)^2\toup{\sigma_1\otimes\sigma_2} \det(M: g_1M)\otimes \det(M: g_2M)
\toup{\id\otimes g_1} \det(M:g_1M)\otimes\det(g_1M: g_1g_2M)\\ 
\iso\, \det(M: g_1g_2M)
\end{multline*}

Informally speaking, one may think of the exact sequence of group stacks 
$$
1\to B(\mu_2)\to \wt G(F)\to G(F)\to 1
$$ 
We also have a canonical section $G(\cO)\to \wt G(F)$ sending $g$ to 
$$
(g,\cB=k, \id:\cB^2\,\iso\, \det(M:M))
$$ 

 The group stack $\wt G(F)$ acts naturally on $\wt\cL_d(M(F))$, the action map $\wt G(F)\times\wt\cL_d(M(F))\to \wt\cL_d(M(F))$ sends 
$$
(g, \cB_1, \sigma_1: \cB_1^2\,\iso\, \det(M:gM)), (L, \cB_2, \sigma_2: \cB_2^2\,\iso\, \det(M:L))
$$
to the collection $(gL, \cB)$, where $\cB=\cB_1\otimes\cB_2$ is equipped with the isomorphism
$$
(\cB_1\otimes\cB_2)^2\toup{\sigma_1\otimes\sigma_2} \det(M:gM)\otimes \det(M:L)\toup{\id\otimes g}\det(M:gM)\otimes\det(gM:gL)\,\iso\, \det(M:gL)
$$

\noindent
5.5 For $g\in G(F)$ and a c-lattice $R\subset R^{\perp}$ in $M(F)$ we have an isomorphism of symplectic spaces $g: R^{\perp}/R\,\iso\, (gR)^{\perp}/gR$. For each c-lattice $R_1\subset R$ we have a diagram
$$
\begin{array}{ccc}
\cL(R_1^{\perp}/R_1)_R & \isoup{g} & \cL(gR_1^{\perp}/gR_1)_{gR}\\
\downarrow\lefteqn{\scriptstyle p} && \downarrow\lefteqn{\scriptstyle p}\\
\cL(R^{\perp}/R) & \isoup{g} & \cL(gR^{\perp}/gR)
\end{array}
$$
Let $\cA_{R_1}$ be the ($\ZZ/2\ZZ$-graded purely of degree zero) line bundle on $\cL(R_1^{\perp}/R_1)_R$ whose fibre at $L$ is $\det L\otimes\det(M: R_1)$. Assume that $\tilde g=(g, \cB, \cB^2\,\iso\, \det(M: gM))$ is a $k$-point of $\tilde G(F)$ over $g$. It yields a diagram
$$
\begin{array}{ccc}
\wt\cL(R_1^{\perp}/R_1)_R & \isoup{\tilde g} & \wt\cL(gR_1^{\perp}/gR_1)_{gR}\\
\downarrow\lefteqn{\scriptstyle p} && \downarrow\lefteqn{\scriptstyle p}\\
\wt\cL(R^{\perp}/R) & \isoup{\tilde g} & \wt\cL(gR^{\perp}/gR)
\end{array}
$$
Here the top horizontal arrow sends $(L, \cB_1, \cB_1^2\,\iso\, \det L\otimes\det(M:R_1))$ to 
$$
(gL, \cB_2, \sigma:\cB_2^2\,\iso\, \det(gL)\otimes\det(M: gR_1)),
$$ 
where $\cB_2=\cB_1\otimes\cB$ and $\sigma$ is the composition
\begin{multline*}
(\cB_1\otimes\cB)^2\,\iso\, \det L\otimes\det(M: R_1)\otimes \det(M:gM)\,\toup{g\otimes\id}\\
\det (gL)\otimes\det(gM: gR_1)\otimes \det(M:gM)\,\iso\,
\det (gL)\otimes \det(M:gR_1)
\end{multline*}
In the limit by $R_1$ the corresponding functors $\tilde g^*: \P(\wt\cL(gR_1^{\perp}/gR_1)_{gR})\,\iso\, \P(\wt\cL(R_1^{\perp}/R_1)_R)$ yield an equivalence
$$
\tilde g^*: \P(\wt\cL_d(M(F))_{gR})\,\iso\, \P(\wt\cL_d(M(F))_R)
$$
Taking one more limit by the partially ordered set of c-lattices $R$, one gets an equivalence 
$$
\tilde g^*: \P(\wt\cL_d(M(F)))\,\iso\, \P(\wt\cL_d(M(F)))
$$ 
In this sense $\tilde G(F)$ acts on $\P(\wt\cL_d(M(F)))$.

\bigskip\bigskip

\centerline{\scshape 6. Canonical interwining operators: local field case}

\bigskip\noindent
6.1 Keep notations of Section~5. Write $H=M\oplus\Omega$ for the Heisenberg group defined as in Section~2.1, this is a group scheme over $\Spec\cO$. 

 For $L\in \cL_d(M(F))$ we have the subgroup $\bar L=L\oplus \Omega(F)\subset H(F)$ and the character $\chi_L: \bar L\to\Qlb^*$ given by $\chi_L(l,a)=\chi(a)$. Here $\chi:\Omega(F)\to\Qlb^*$ sends $a$ to $\psi(\Res a)$. In the classical setting we let $\cH_L$ denote the space of functions $f: H(F)\to\Qlb$ satisfying
\begin{itemize} 
\item[C1)] $f(\bar l h)=\chi_L(\bar l)f(h)$, for $h\in H, \bar l\in\bar L$;
\item[C2)] there exists a c-lattice $R\subset M(F)$ such that 
$f(h(r,0))=f(h)$ for $r\in R, h\in H$. 
\end{itemize}
Note that such $f$ has automatically compact support modulo $\bar L$. The group $H(F)$ acts on $\cH_L$ by right translations, this is a model of the Weil representation. Let us introduce a geometric analog of $\cH_L$. 

 Given a c-lattice $R\subset M(F)$ such that $R\subset R^{\perp}$ write $H_R=(R^{\perp}/R)\oplus k$ for the Heisenberg group corresponding to the symplectic space $R^{\perp}/R$. If $L\in \cL_d(M(F))_R$ then $L_R:=L\cap R^{\perp}\subset R^{\perp}/R$ is lagrangian. Set $\bar L_R=L_R\oplus k\subset H_R$. Let $\chi_{L,R}: \bar L_R\to \Qlb^*$ be the character sending $(l,a)$ to $\psi(a)$. Set 
$$
\cH_{L_R}=\{f: H_R\to \Qlb\mid f(\bar l h)=\chi_{L,R}(\bar l)f(h),\; h\in H_R, \bar l\in \bar L_R\}
$$  

\begin{Lm} There is a canonical embedding $T^L_R: \cH_{L_R}\hook{}\cH_L$ whose image is the subspace of those $f\in\cH_L$ which satisfy 
\begin{equation}
\label{cond_equiv_R_appendix}
f(h(r,0))=f(h)\;\; \mbox{for} \;\;r\in R, h\in H
\end{equation}
\end{Lm}
\begin{Prf}
Set 
$$
'\cH_{L_R}=\{\phi: R^{\perp}/R\to\Qlb\mid \phi(r+l)=\chi(\frac{1}{2}\omega\<r,l\>)\phi(r), \; r\in R^{\perp}/R, \; l\in L_R\}
$$
We have an isomorphism $\cH_{L_R}\,\iso\, {'\cH_{L_R}}$ sending $f$ to $\phi$ given by $\phi(r)=f(r,0)$. Given $f\in \cH_L$ satisfying (\ref{cond_equiv_R_appendix}), we associate to $f$ a function $\phi\in{'\cH_{L_R}}$ given by 
$$
\phi(r)=q^{\frac{1}{2}\dim R^{\perp}/R}
f(r,0)
$$ 
for $r\in R^{\perp}$. This defines the map $T^L_R$.
\end{Prf}

\medskip

 Assume that $S\subset R\subset M(F)$ are c-lattices and $R\cap L=0$. Remind the operator $\cH_{L_R}\toup{T^{L_S}} \cH_{L_S}$ given by (\ref{def_T^L_classical}), it corresponds to the isotropic subspace $R/S\subset S^{\perp}/S$. The composition 
$\cH_{L_R}\toup{T^{L_S}} \cH_{L_S}\toup{T^L_S} \cH_L$ equals $T^L_R$. 

 The geometric analog of $\cH_L$ is as follows. For a c-lattice $R$ such that $R\cap L=0$ and $R\subset R^{\perp}$ we have the category $\cH_{L_R}$ of perverse sheaves on $H_R$ which are $(\bar L_R, \chi_{L,R})$-equivariant, and the corresponding category $\D\!\cH_{L_R}$. For $S\subset R$ as above we have an (exact for the perverse structure and fully faithful) transition functor (\ref{def_T^L_geom}), which we now denote by 
$$
T^L_{S,R}: \D\!\cH_{L_R}\to \D\!\cH_{L_S}
$$
Define $\cH_L$ (resp., $\D\!\cH_L$) as the inductive 2-limit of $\cH_{L_R}$ (resp., of $\D\!\cH_{L_R}$) over the partially ordered set of c-lattices $R$ such that $R\cap L=0$ and $R\subset R^{\perp}$. So, $\cH_L$ is abelian and $\D\!\cH_L$ is a triangulated category. 

\medskip\noindent
6.2 Let $R\subset R^{\perp}$ be a c-lattice in $M(F)$. We have a projection 
$$
\cL_d(M(F))_R\to \cL(R^{\perp}/R)
$$ 
sending $L$ to $L_R$.  Let $\cA_R$ be the $\ZZ/2\ZZ$-graded purely of degree zero line bundle on $\cL(R^{\perp}/R)$ whose fibre at $L_1$ is $\det L_1\otimes\det(M:R)$. Write $\wt\cL(R^{\perp}/R)$ for the gerb of square roots of $\cA_R$. The restriction of $\cA_R$ to $\cL_d(M(F))_R$ identifies canonically with $\cA_d$. 
The above projection lifts naturally to a morphism of gerbs
\begin{equation}
\label{map_gerbs_R_appendix}
\wt\cL_d(M(F))_R\to \wt\cL(R^{\perp}/R)
\end{equation} 

Given $k$-points $N^0, L^0\in \wt\cL_d(M(F))$ we are going to associate to them in a canonical way a functor 
\begin{equation}
\label{functor_cF_local_field_case}
\cF_{N^0,L^0}: \D\!\cH_L\to\D\!\cH_N
\end{equation}
sending $\cH_L$ to $\cH_N$. To do so, consider a c-lattice $R\subset R^{\perp}$ in $M(F)$ such that $L,N\in \cL_d(M(F))_R$. Write $N_R^0, L_R^0\in {\wt\cL(R^{\perp}/R)}$ for the images of $N^0$ and $L^0$ under (\ref{map_gerbs_R_appendix}). By definition, the enhanced structure on $L_R$ and $N_R$ is given by one-dimensional vector spaces $\cB_L, \cB_N$ equipped with 
$$
\cB_L^2\,\iso\, \det L_R\otimes\det(M:R),\;\;\; 
\cB_N^2\,\iso\, \det N_R\otimes\det(M:R),
$$
hence an isomorphism $\cB^2\,\iso\, \det L_R\otimes\det N_R$ for $\cB:=\cB_L\otimes\cB_N\otimes\det(M:R)^{-1}$. 
We denote by 
$$
\cF_{N_R^0, L_R^0}: \D\!\cH_{L_R}\to \D\!\cH_{N_R}
$$ 
the canonical interwining functor defined in Section~3.5 corresponding to $(N_R, L_R,\cB)\in\tilde Y$, here $Y=\cL(R^{\perp}/R)\times\cL(R^{\perp}/R)$. The following is an immediate consequence of Proposition~\ref{Pp_appendix_compatibility_property_CIO}. 

\begin{Pp}  Let $S\subset R\subset R^{\perp}\subset S^{\perp}$ be c-lattices such that $L^0,N^0\in \wt\cL_d(M(F))_R$. Then the following diagram of categories is canonically 2-commutative
$$
\begin{array}{ccccccc}
\D\!\cH_{L_R} &&& \toup{T^L_{S,R}} &&& \D\!\cH_{L_S}\\
\downarrow\lefteqn{\scriptstyle \cF_{N_R^0, L_R^0}} &&&&&& \downarrow\lefteqn{\scriptstyle \cF_{N_S^0, L_S^0}}\\
\D\!\cH_{N_R} &&& \toup{T^N_{S,R}} &&& \D\!\cH_{N_S}
\end{array}
$$
\end{Pp}

 Define (\ref{functor_cF_local_field_case}) as the limit of functors $\cF_{N^0_R, L^0_R}$ over the partially ordered set of c-lattices $R\subset R^{\perp}$ such that $L,N\in \cL_d(M(F))_R$. As in Section~3.5, one shows that for $L^0,N^0,R^0\in \wt\cL_d(M(F))$ the diagram is canonically 2-commutative
$$
\begin{array}{ccc}
\D\!\cH_L & \;\toup{\cF_{R^0,L^0}} & \;\;\;\;\D\!\cH_R\\
 & \searrow\lefteqn{\scriptstyle \cF_{N^0, L^0}} & \;\;\;\;\downarrow\lefteqn{\scriptstyle \cF_{N^0, R^0}}\\
&& \;\;\;\;\D\!\cH_N
\end{array}
$$
 
 Our main result in the local field case is as follows.

\begin{Th} 
\label{Th_2}
For each $k$-point $L^0\in \wt\cL_d(M(F))$ there is a canonical functor 
\begin{equation}
\label{functor_cF_L^0_for_Th}
\cF_{L^0}: \D\!\cH_L\to \D(\wt\cL_d(M(F)))
\end{equation}
sending $\cH_L$ to $\P(\wt\cL_d(M(F)))$. For a pair of $k$-points $(L^0,N^0)$ in $\wt\cL_d(M(F))$ the diagram
\begin{equation}
\label{diag_compatibility_functors_cF_local_field}
\begin{array}{ccc}
\D\!\cH_L & \toup{\cF_{L^0}} & \D(\wt\cL_d(M(F)))\\
\downarrow\lefteqn{\scriptstyle \cF_{N^0,L^0}} &\;\;\;\;\;\;\;\;\nearrow\lefteqn{\scriptstyle \cF_{N^0}}\\
\D\!\cH_N
\end{array}
\end{equation}
is canonically 2-commutative. Let $W(\wt\cL_d(M(F)))$
be the essential image of 
$$
\cF_{L^0}: \cH_L\to \P(\wt\cL_d(M(F))),
$$
this is a full subcategory independent of $L^0$. Besides, $W(\wt\cL_d(M(F)))$ is preserved under the natural action of $\tilde G(F)$ on $\P(\wt\cL_d(M(F)))$. 
\end{Th}

 We will refer to $W(\wt\cL_d(M(F))$ as \select{the non-ramified Weil category on} $\wt\cL_d(M(F))$. Remind that in the classical setting 
$$
\cH_L=\cH_{L,odd}\oplus \cH_{L,even}
$$ 
is a direct sum of two irreducible representations of the metaplectic group (consisting of odd and even functions respectively). The representation $\cH_{L,odd}$ is ramified, whence $\cH_{L,even}$ is not. The category $W(\wt\cL_d(M(F)))$ together with the action of $\tilde G(F)$ is a geometric counterpart of the representation $\cH_{L,even}$. The proof of Theorem~\ref{Th_2} is given in Sections~6.3-6.4.

\medskip\noindent
6.3 Let $L^0$ be a $k$-point of $\wt\cL_d(M(F))$. 
Let $R\subset R^{\perp}$ be a c-lattice with $L\cap R=0$. Write $L^0_R$ for the image of $L^0$ under (\ref{map_gerbs_R_appendix}). Applying the construction of Section~3.6 to the symplectic space $R^{\perp}/R$ with $L_R^0\in\tilde\cL(R^{\perp}/R)$, one gets the functor 
$$
\cF_{L_R^0}: \D\!\cH_{L_R}\to \D(\wt\cL(R^{\perp}/R))
$$ 

 If $N^0$ is another $k$-point of $\wt\cL_d(M(F))_R$ then writing $N^0_R$ for the image of $N^0$ in $\wt\cL(R^{\perp}/R)$ we also get that the diagram
\begin{equation}
\label{diag_one_more_commutativity_CIO}
\begin{array}{ccc}
\D\!\cH_{L_R} & \toup{\cF_{L_R^0}} & \D(\wt\cL(R^{\perp}/R))\\
\downarrow\lefteqn{\scriptstyle \cF_{N^0_R, L^0_R}} & \;\;\;\;\;\;\;\nearrow\lefteqn{\scriptstyle \cF_{N_R^0}}\\
\D\!\cH_{N_R} 
\end{array}
\end{equation}
is canonically 2-commutative. 

 Let now
$$
_R\cF_{L^0}: \D\!\cH_{L_R}\to \D(\wt\cL_d(M(F))_R)
$$
denote the composition of $\cF_{L_R^0}$ with the (exact for the perverse t-structures) restriction functor $\D(\wt\cL(R^{\perp}/R))\to \D(\wt\cL_d(M(F))_R)$ for the projection  (\ref{map_gerbs_R_appendix}). 

 Let $S\subset R$ be another c-lattice. As in Section~5.3, for the open immersion $j_{S,R}: \wt\cL_d(M(F))_R\hook{} \wt\cL_d(M(F))_S$ we have the restriction functors $j^*_{S,R}: \D(\wt\cL_d(M(F))_S)\to \D(\wt\cL_d(M(F))_R)$. 
 
\begin{Lm} 
\label{Lm_comaptibility_restrictions_appendix}
The diagram of functors is canonically 2-commutative
$$
\begin{array}{ccc}
\D\!\cH_{L_R} & \toup{_R\cF_{L^0}} & \D(\wt\cL_d(M(F))_R)\\
\downarrow\lefteqn{\scriptstyle  T^L_{S,R}} && \uparrow\lefteqn{\scriptstyle j_{S,R}^*}\\
\D\!\cH_{L_S} & \toup{_S\cF_{L^0}} & \D(\wt\cL_d(M(F))_S)
\end{array}
$$
\end{Lm} 
\begin{Prf}
We have an open immersion $j: \wt\cL(S^{\perp}/S)_R\hook{} \wt\cL(S^{\perp}/S)$ and a projection $p_{R/S}: \wt\cL(S^{\perp}/S)_R)\to \wt\cL(R^{\perp}/R)$. Set $P_{R/S}=p_{R/S}^*\otimes(\Qlb[1](\frac{1}{2}))^{\dimrel(p_{R/S})}$. 
It suffices to show that the diagram is canonically 2-commutative
$$
\begin{array}{ccccc}
\D\!\cH_{L_R} & \;\toup{\cF_{L_R^0}} & \D(\wt\cL(R^{\perp}/R)) & \toup{P_{R/S}} &
\D(\wt\cL(S^{\perp}/S)_R))\\
\downarrow\lefteqn{\scriptstyle  T^L_{S,R}} &&&
\nearrow\lefteqn{\scriptstyle j^*}\\
\D\!\cH_{L_S} & \;\toup{\cF_{L^0_S}} & \D(\wt\cL(S^{\perp}/S))  
\end{array}
$$
This follows from Lemma~\ref{Lm_appendix_compatibility_property}. 
\end{Prf}
 
\medskip

 Define $\cF_{L^0, R}: \D\!\cH_{L_R}\to \D(\wt\cL_d(M(F)))$ as the functor sending $K_1$ to the following object $K_2$. For a c-lattice $S\subset R$ we declare the restriction of $K_2$ to
$\wt\cL_d(M(F))_S$ to be 
$$
(_S\cF_{L^0} \comp T^L_{S,R})(K_1)
$$ 
By Lemma~\ref{Lm_comaptibility_restrictions_appendix}, the corresponding projective system defines an object $K_2$ of $\D(\wt\cL_d(M(F)))$. 
 
 Finally, for $S\subset R$ with $R\cap L=0$ 
the diagram
$$
\begin{array}{ccc}
\D\!\cH_{L_R} & \toup{\cF_{L^0, R}} & \D(\wt\cL_d(M(F)))\\
\downarrow\lefteqn{\scriptstyle T^L_{S,R}} & \;\;\;\nearrow\lefteqn{\scriptstyle \cF_{L^0, S}}\\
\D\!\cH_{L_S} 
\end{array}
$$
is canonically 2-commutative. 
We define (\ref{functor_cF_L^0_for_Th}) as the limit of the functors $\cF_{L^0,R}$ over the partially ordered set of c-lattices $R\subset R^{\perp}$ such that $L\cap R=0$. 
The commutativity of (\ref{diag_compatibility_functors_cF_local_field}) follows from the commutativity of (\ref{diag_one_more_commutativity_CIO}).  
 
\begin{Def} The non-ramified Weil category $W(\wt\cL_d(M(F)))$ is the essential image of the functor $\cF_{L^0}: \cH_L\to \P(\wt\cL_d(M(F)))$. It does not depend on a choice of a $k$-point $L^0$ of $\wt\cL_d(M(F))$. 
\end{Def} 
 
\noindent
6.4 Let $R\subset R^{\perp}$ be a c-lattice in $M(F)$, let $\tilde g\in\tilde G(F)$ be a $k$-point, write $g$ for its image in $G(F)$. As in Section~5.5, we have an isomorphism $g: H_R\,\iso\, H_{gR}$ of algebraic groups over $k$ sending $(x,a)\in (R^{\perp}/R)\times\A^1$ to $(gx,a)\in (gR^{\perp}/gR)\times\A^1$. For $L\in \cL_d(M(F))_R$ it induces an equivalence
$$
g: \cH_{L_R}\,\iso\, \cH_{gL_{gR}}
$$
If $L^0\in\wt\cL_d(M(F))_R$ is a $k$-point then the $G$-equivariance of $F$ implies that the diagram is canonically 2-commutative
$$
\begin{array}{ccc}
\cH_{L_R} & \toup{\cF_{L^0_R}} & \P(\wt\cL(R^{\perp}/R))\\
\downarrow\lefteqn{\scriptstyle g} && \downarrow\lefteqn{\scriptstyle \tilde g} \\
\cH_{gL_{gR}} & \toup{\cF_{\tilde gL^0_{gR}}} & \P(\wt\cL(gR^{\perp}/gR))
\end{array}
$$
This, in turn, implies that the diagram is 2-commutative
$$
\begin{array}{ccc}
\cH_{L_R} & \toup{\cF_{L^0,R}} & \P(\wt\cL_d(M(F)))\\
\downarrow\lefteqn{\scriptstyle g} && \downarrow\lefteqn{\scriptstyle \tilde g}\\
\cH_{gL_{gR}} & \toup{\cF_{\tilde gL^0, gR}} & \P(\wt\cL_d(M(F)))
\end{array}
$$
Thus, Theorem~\ref{Th_2} is proved. 
 
\medskip\noindent
6.5 {\scshape Theta-sheaf\ } Let $L\in \cL_d(M(F))_M$, this is equivalent to saying that $L\subset M(F)$ is a lagrangian d-lattice such that $L\oplus M=M(F)$. Then the category $\cH_{L_M}$ has a distinguished object $\cL_{\psi}$ on $\A^1=\H_M$. Write $S_L$ for its image under $\cH_{L_M}\to\cH_L$. The line bundle $\cA_d$ over $\cL_d(M(F))_M$ is canonically trivialized, so $L$ has a dintinguished enhanced structure 
$$
(L,\cB)=L^0\in \wt\cL_d(M(F))_M,
$$
 where $\cB=k$ is equipped with $\id: \cB^2\,\iso\,\det(M:L)$. The \select{theta-sheaf $S_{M(F)}$ over} $\wt\cL_d(M(F))$ is defined as $\cF_{L^0}(S_L)$. It does not depend on $L\in \cL_d(M(F))_M$ in the sense that for another $N\in  \cL_d(M(F))_M$ the diagram (\ref{diag_compatibility_functors_cF_local_field}) yields a canonical isomorphism $\cF_{L^0}(S_L)\,\iso\, \cF_{N^0}(S_N)$. The perverse sheaf $S_{M(F)}$ has a natural $G(\cO)$-equivariant structure. 
 
\medskip\noindent
6.6 {\scshape Relation with the Schr\" odinger model}

\medskip\noindent
Assume in addition that $M$ is decomposed as $M\,\iso\, U\oplus U^*\otimes\Omega$, where $U$ is a free $\cO$-module of rank $d$, both $U$ and $U^*\otimes\Omega$ are isotropic, and the form $\omega: \wedge^2 M\to\Omega$ is given by $\omega\<u,u^*\>=\<u, u^*\>$ for $u\in U, u^*\in U^*\otimes\Omega$, where $\<\cdot, \cdot\>$ is the natural pairing between $U$ and $U^*$. 
Let $\bar U=U(F)\oplus\Omega(F)$ viewed as a subgroup of $H(F)$, it is equipped with the character $\chi_U: \bar U\to \Qlb^*$ given by $\chi_U(u,a)=\psi(\Res a)$, $a\in\Omega(F), u\in U(F)$. Write 
\begin{multline*}
\Shr_U=\{f: H(F)\to\Qlb\mid f(\bar uh)=\chi_U(\bar u)f(h), \bar u\in \bar U, h\in H(F),\;\; f\;\;\mbox{is smooth}, \\ \mbox{of compact support modulo}\; \bar U 
\},
\end{multline*}
$H(F)$ acts on it by right translations. This is the Schr\" odinger model of the Weil representation, it identifies naturally with the Schwarz space 
$\cS(U^*\otimes\Omega(F))$.  

 Remind the definition of the derived category $\D(U^*\otimes\Omega)$ and its subcategory of perverse sheaves $\P(U^*\otimes\Omega)$ given in (\cite{L2}, Section~4). For $N,r\in \ZZ$ with $N+r\ge 0$ we write $_{N,r}U=t^{-N}U/t^r U$. 
 
 For $N_1\ge N_2, r_1\ge r_2$ we have a diagram 
$$
_{N_2, r_2}(U^*\otimes\Omega) \getsup{p} {_{N_2, r_1}(U^*\otimes\Omega)}\toup{i} {_{N_1,r_1}(U^*\otimes\Omega)},
$$
 where $p$ is the smooth projection and $i$ is a closed immersion. We have a transition functor 
\begin{equation}
\label{functor_transition_for_U}
\D(_{N_2,r_2}(U^*\otimes\Omega))\to \D(_{N_1,r_1}(U^*\otimes\Omega))
\end{equation}
sending $K$ to $i_! p^*K\otimes(\Qlb[1](\frac{1}{2}))^{\dimrel(p)}$, 
it is fully faithful and exact for the perverse t-structures. Then $\D(U^*\otimes\Omega(F))$ (resp., $\P(U^*\otimes\Omega(F))$) is defined as the inductive 2-limit of $\D(_{N,r}(U^*\otimes\Omega))$ (resp., of $\P(_{N,r}(U^*\otimes\Omega))$) as $r,N$ go to infinity. 
The category $\P(U^*\otimes\Omega(F))$ is the geometric analog of the space $\Shr_U$. 
 
 In this section we prove the following. 
 
\begin{Pp}  
\label{Pp_for_Shrodinger}
For each $k$-point $L^0\in\wt\cL_d(M(F))$ there is a canonical equivalence
\begin{equation}
\label{functor_cF_for_U(F)}  
\cF_{U(F), L^0}: \D(U^*\otimes\Omega(F))\to \D\!\cH_L
\end{equation}  
which identifies $\P(U^*\otimes\Omega(F))$ with the category $\cH_L$. For $L^0, N^0\in \wt\cL_d(M(F))$ the diagram is canonically 2-commutative
$$
\begin{array}{ccc}
\D(U^*\otimes\Omega(F)) &\toup{\cF_{U(F), L^0}} & \D\!\cH_L \\
\downarrow\lefteqn{\scriptstyle \cF_{U(F), N^0}} & \nearrow\lefteqn{\scriptstyle \cF_{L^0, N^0}}\\
\D\!\cH_N
\end{array}
$$
\end{Pp}

 For $N\ge 0$ consider the c-lattice $R=t^NM$ in $M(F)$ and the corresponding symplectic space $R^{\perp}/R={_{N,N}M}$. Set $U_R:={_{N,N}U}\in \cL(_{N,N}M)$. We have the line bundle $\cA_N$ on $\cL(_{N,N}M)$ whose fibre at $L$ is $\det(_{0,N}M)\otimes\det L$. As above, $\wt\cL(_{N,N}M)$ is the gerb of square roots of $\cA_N$. Let 
$$
U^0_R=(U_R, \det(_{0,N}U))\in \wt\cL(_{N,N}M)
$$
equipped with a canonical $\ZZ/2\ZZ$-graded isomorphism 
$\det(_{0,N}U)^2\,\iso\, \det U_R\otimes \det(_{0,N}M)$.

  Let $H_R={_{N,N}M}\times\A^1$ denote the corresponding Heisenberg group, it has the subgroup $\bar U_R=U_R\times\A^1$ equipped with the character
$\chi_{U,R}: \bar U_R\to\Qlb^*$ given by $\chi_{U,R}(u,a)=\psi(a)$, $a\in\A^1$. In the classical setting, $\cH_{U_R}$ is the space of functions on $H_R$, which are $(\bar U_R, \chi_{U,R})$-equivariant under the left multiplication. 
Set $\Shr_U^R=\{f\in\Shr_U\mid f(h(r,0))=f(h), r\in R, h\in H\}$. 
  
\begin{Lm}  In the classical setting there is an isomorphism 
\begin{equation}
\label{iso_U_level_N}
\Shr_U^R\,\iso\, \cH_{U_R}
\end{equation} 
\end{Lm} 
\begin{Prf}  Write $\cH'_{U_R}=\{\phi': R^{\perp}/R\to\Qlb\mid \phi'(m+u)=\psi(\frac{1}{2}\<m,u\>)\phi'(m), \; u\in U_R \}$. We identify $\cH_{U_R}\,\iso\,\cH'_{U_R}$ via the map $\phi\mapsto \phi'$, where $\phi'(m)=\phi(m,0)$. Given $f\in \Shr_U^R$ for $m\in t^{-N}M$ the value $f(m,0)$ depends only on the image $\bar m$ of $m$ under $t^{-N}M\to {_{N,N}M}$. The isomorphism (\ref{iso_U_level_N}) sends $f$ to $\phi'\in \cH'_{U_R}$ given by $\phi(\bar m)=f(m,0)$.
\end{Prf}
 
\medskip

 In the geometric setting $\cH_{U_R}$ is the category of  $(\bar U_R, \chi_{U,R})$-equivariant perverse sheaves on $H_R$. We identify it with $\P(_{N,N}(U^*\otimes\Omega))$ as follows. Let $m_U: 
\bar U_R\times {_{N,N}(U^*\otimes\Omega)}\to H_R$ be the isomorphism sending $(\bar u, h)$ to their product $\bar u h$ in $H_R$. The functor $\D(_{N,N}(U^*\otimes\Omega))\to \D\!\cH_{U_R}$ sending $K$ to 
$$
(m_U)_!(\chi_{U,R}\boxtimes K)\otimes(\Qlb[1](\frac{1}{2}))^{\dim \bar U_R}
$$
is an equivalence (exact for the perverse t-structures).
 
 Let $N'\ge N$ and $S=t^{N'}M$. The corresponding transition functor (\ref{functor_transition_for_U}) now yields a functor denoted $T^U_{S,R}: \D\!\cH_{U_R}\to\D\!\cH_{U_{S}}$. 
 
 Let $L^0\in \wt\cL_d(M(F))$ be a $k$-point over $L\in \cL_d(M(F))$. Assume that $N$ is large enough so that $L\cap R=0$. Let $L_R^0$ denote the image of $L^0$ under (\ref{map_gerbs_R_appendix}). Define $U_S^0, L^0_S\in\wt\cL(S^{\perp}/S)$ similarly.
 
\begin{Lm} 
\label{Lm_for_Shrodinger}
The diagram is canonically 2-commutative
$$
\begin{array}{ccc} 
\D\!\cH_{U_R}\;\;\;\;\; & \,\toup{T^U_{S,R}} & \;\;\;\;\D\!\cH_{U_{S}}\\
\downarrow\lefteqn{\scriptstyle \cF_{L^0_R, U^0_R}} \;\;\;\;\;&& \;\;\;\; \downarrow\lefteqn{\scriptstyle \cF_{L^0_{S}, U^0_{S}}}\\
\D\!\cH_{L_R}\;\;\;\;\; & \,\toup{T^L_{S,R}} & \;\;\;\;\D\!\cH_{L_{S}}
\end{array}
$$
\end{Lm}
\begin{Prf}  
Set $W=t^{N'}U\oplus t^N(U^*\otimes\Omega)$. The subspace $W/S\subset S^{\perp}/S$ is isotropic, and $U_S\cap (W/S)=L_S\cap (W/S)=0$. Write $H_W=(W^{\perp}/W)\times \A^1$ for the corresponding Heisenberg group. Set $U_W=U_S\cap (W^{\perp}/S)$, $L_W=L_S\cap (W^{\perp}/S)$. 
Applying Proposition~\ref{Pp_appendix_compatibility_property_CIO}, we get a 2-commutative diagram
$$
\begin{array}{ccc}
\D\!\cH_{U_W}\;\;\;\;\; & \,\toup{T^U_{S,W}} & \;\;\;\;\D\!\cH_{U_S}\\
\downarrow\lefteqn{\scriptstyle \cF_{L^0_W, U^0_W}} \;\;\;\;\;&& \;\;\;\; \downarrow\lefteqn{\scriptstyle \cF_{L^0_{S}, U^0_{S}}}\\
\D\!\cH_{L_W}\;\;\;\;\; & \,\toup{T^L_{S,W}} & \;\;\;\;\D\!\cH_{L_{S}}
\end{array}
$$
Now $R/W\subset W^{\perp}/W$ is an isotropic subspace, and $R/W\subset U_W$, $R/W\cap L_W=0$. Note that $U_R=U_W/(R/W)$. Applying  Proposition~\ref{Pp_second_compatibility}, we get a 2-commutative diagram
$$
\begin{array}{ccc} 
\D\!\cH_{U_R}\;\;\;\;\; & \,\toup{T^U_{W,R}} & \;\;\;\;\D\!\cH_{U_{W}}\\
\downarrow\lefteqn{\scriptstyle \cF_{L^0_R, U^0_R}} \;\;\;\;\;&& \;\;\;\; \downarrow\lefteqn{\scriptstyle \cF_{L^0_{W}, U^0_{W}}}\\
\D\!\cH_{L_R}\;\;\;\;\; & \,\toup{T^L_{W,R}} & \;\;\;\;\D\!\cH_{L_{W}}
\end{array}
$$
Our assertion easily follows.
\end{Prf} 
 
\medskip

\begin{Prf}\select{of Proposition~\ref{Pp_for_Shrodinger}}

\smallskip\noindent
Passing to the limit as $N$ goes to infinity, the functors $\cF_{L^0_R, U^0_R}: \D\!\cH_{U_R}\to \D\!\cH_{L_R}$ from Lemma~\ref{Lm_for_Shrodinger} yield the desired functor (\ref{functor_cF_for_U(F)}). The second assertion follows by construction. 
\end{Prf}
   
\begin{Def} Let $\cF_{U(F)}: \D(U^*\otimes\Omega(F))\to \D(\wt\cL_d(M(F)))$ denote the composition  
$$
\D(U^*\otimes\Omega(F))
\toup{\cF_{U(F), L^0}}\D\!\cH_L \toup{\cF_{L^0}} \D(\wt\cL_d(M(F)))
$$ 
By Theorem~\ref{Th_2} and Proposition~\ref{Pp_for_Shrodinger}, it does not depend on the choice of a $k$-point $L^0\in \wt\cL_d(M(F))$. By construction, $\cF_{U(F)}$ is exact for the perverse t-structures. 
\end{Def}  
   
   We have a morphism of group stacks $\GL(U)(F)\to \tilde G(F)$ sending $g\in\GL(U)(F)$ to $(g, \cB=\det(U: gU))$ equipped with a canonical $\ZZ/2\ZZ$-graded isomorphism
$$
\det(M: gM)\,\iso\,\det(U: gU)\otimes \det(U^*\otimes\Omega: g(U^*\otimes\Omega))\,\iso\, \det(U: gU)^{\otimes 2}
$$    
Let $\GL(U)(F)$ act on $\wt\cL_d(M(F))$ via this homomorphism, let it also act naturally on $U^*\otimes\Omega(F)$. Then one may show that $\cF_{U(F)}$ commutes with the action of $\GL(U)(F)$. 
  
  Note also that over $\GL(U)(\cO)$ the sections $\GL(U)(F)\to \tilde G(F)$ and $G(\cO)\to \tilde G(F)$ are compatible.
\bigskip\medskip

\centerline{\scshape 7. Global application}

\bigskip\noindent
7.1  Assume $k$ algebraically closed. Let $X$ be a smooth connected projective curve. Let $\Omega$ be the canonical invertible sheaf on $X$. Let $G$ be the group scheme over $X$ of automorphisms of $\cO_X^d\oplus \Omega^d$ perserving the symplectic form $\wedge^2(\cO_X^d\oplus \Omega^d)\to\Omega$. 

 Write $\Bun_G$ for the stack of $G$-torsors on $X$, it classifies a rank $2d$-vector bundle $\cM$ on $X$ together with a symplectic form $\wedge^2 \cM\to\Omega$. Let $\cA$ be the ($\ZZ/2\ZZ$-graded purely of degree zero) line bundle on $\Bun_G$ whose fibre at $\cM$ is $\det\RG(X,\cM)$. Write $\Bunt_G$ for the gerb of square roots of $\cA$ over $\Bun_G$. 
 
 Remind the definition of the theta-sheaf $\Aut$ on $\Bunt_G$ (\cite{L1}, Definition~1). Let $_i\Bun_G\hook{}\Bun_G$ be the locally closed substack given by $\dim\H^0(X,\cM)=i$ for $\cM\in \Bun_G$. Write $_i\Bunt_G$ for the restriction of $\Bunt_G$ to $_i\Bun_G$. 
 
 Let $_i\cB$ be the line bundle on $_i\Bun_G$ whose fibre at $\cM\in {_i\Bun_G}$ is $\det\H^0(X,\cM)$, we view it as $\ZZ/2\ZZ$-graded of degree $i\!\mod 2$. For each $i$ we have a canonical $\ZZ/2\ZZ$-graded isomorphism
$_i\cB^2\,\iso\,\cA$, it yields a trivialization $_i\Bunt_G\,\iso\, {_i\Bun_G}\times B(\mu_2)$.   
 
 Define $\Aut_g\in \P(\Bunt_G$ (resp., $\Aut_s\in\P(\Bunt_G)$) as the intermediate extension of 
$$
(\Qlb\boxtimes W)\otimes(\Qlb[1](\frac{1}{2})^{\dim\Bun_G}
$$ 
(resp., of
$(\Qlb\boxtimes W)\otimes (\Qlb[1](\frac{1}{2})^{\dim\Bun_G-1}$) under $_i\Bunt_G\hook{}\Bunt_G$. Set $\Aut=\Aut_g\oplus\Aut_s$. 
 
\medskip\noindent 
7.2 Fix a closed point $x\in X$. Write $\cO_x$ for the completed local ring of $X$ at $x$, $F_x$ for its fraction field. Fix a $G$-torsor over $\Spec\cO_x$, we think of it as a free $\cO_x$-module $M$ of rank $2d$ with symplectic form $\wedge^2 M\to\Omega(\cO_x)$ and an action of $G(\cO_x)$. We have a map
$$
\xi_x: \Bun_G\to \cL_d(M(F_x))/G(\cO_x),
$$
where $\cL_d(M(F_x))/G(\cO_x)$ is the stack quotient. It sends $\cM\in\Bun_G$ to the Tate space $\cM(F_x)$ with lagrangian c-lattice $\cM(\cO_x)$ and lagrangian d-lattice $\H^0(X-x, \cM)$. 

 The line bundle $\cA_d$ on $\cL_d(M(F_x))/G(\cO_x)$ is that of Section~5.3. Write $\wt\cL_d(M(F_x))/G(\cO_x)$ for the gerb of square roots of $\cA_d$.  
 
  We have canonically $\xi_x^*\cA_d\,\iso\, \cA$, so $\xi$ lifts naturally to a map of gerbs
$$  
\tilde\xi_x: \Bunt_G\to \wt\cL_d(M(F_x))/G(\cO_x)
$$

 For $r\ge 0$ let $_{rx}\Bun_G\subset\Bun_G$ be the open substack given by $\H^0(X, \cM(-rx))=0$. Write $_{rx}\Bunt_G$ for the restriction of the gerb $\Bunt_G$ to $_{rx}\Bun_G$. If $r'\ge r$ then $_{rx}\Bunt_G\subset {_{r'x}\Bunt_G}$ is an open substack, so we consider the projective 2-limit 
$$
\twolim_{r\to\infty} \D(_{rx}\Bunt_G)
$$  
Note that $\twolim_{r\to\infty} \P(_{rx}\Bunt_G)\,\iso\, \P(\Bunt_G)$ is a full subcategory in the above limit. 
Let us define the restriction functor
\begin{equation}
\label{functor_tilde_xi_x}
\tilde\xi_x^*: \D_{G(\cO)}(\tilde\cL_d(M(F)))\to \twolim_{r\to\infty} \D(_{rx}\Bunt_G)
\end{equation}

To do so, for $N\ge r\ge 0$ and $r_1\ge 2N$ let
\begin{equation}
\label{map_xi_N}
\xi_N: {_{rx}\Bun_G}\to {_r\cL(_{N,N}M)/G(\cO/t^{r_1})}
\end{equation}
be the map sending $\cM$ to the lagragian subspace $\H^0(X, \cM(Nx))\subset {_{N,N}\cM}$. 
If $N_1\ge N\ge r$ and $r_1\ge 2N_1$ then the diagram commutes
$$
\begin{array}{ccc}
_{rx}\Bun_G & \toup{\xi_N} & \;{_r\cL(_{N,N}M)/G(\cO/t^{r_1})}\\
 & \searrow\lefteqn{\scriptstyle \xi_{N_1}} & \;\uparrow\lefteqn{\scriptstyle p}\\
&& \;{_r\cL(_{N_1,N_1}M)/G(\cO/t^{r_1})}
\end{array}
$$
It induces a similar diagram between the gerbs (cf. Section~5.3 for their definition)
$$
\begin{array}{ccc}
_{rx}\Bunt_G & \toup{\tilde\xi_N} & \;(_r\cL(_{N,N}M)/G(\cO/t^{r_1}))\tilde{}\\
 & \searrow\lefteqn{\scriptstyle \tilde\xi_{N_1}} & \;\uparrow\\
&& \;(_r\cL(_{N_1,N_1}M)/G(\cO/t^{r_1}))\tilde{}
\end{array}
$$
The functors $K\mapsto \tilde\xi_N^*K\otimes(\Qlb[1](\frac{1}{2}))^{\dimrel(\xi_N)}$ from $\D_{G(\cO)}(_r\tilde\cL(_{N,N}M))$ to $\D(_{rx}\Bunt_G)$ are compatible with the transition functors, so yield a functor
$$
_r\xi_x^*: \D_{G(\cO)}(_r\tilde\cL_d(M(F)))\to \D(_{rx}\Bunt_G)
$$ 
Passing to the limit by $r$, one gets the desired functor (\ref{functor_tilde_xi_x}). 

\begin{Th} 
\label{Th_3}
The object $\tilde\xi_x^*S_{M(F_x)}$ lies in $\P(\Bunt_G)$, and there is an isomorphism of perverse sheaves
$$
\tilde\xi_x^*S_{M(F_x)}\,\iso\, \Aut
$$
\end{Th}
\begin{Prf}  For $r\ge 0$ consider the map 
$$
\tilde\xi_r: {_{rx}\Bunt_G}\to (\cL(_{r,r}M)/G(\cO/t^{2r}))\tilde{}
$$ 
Set $Y=\cL(_{r,r}M)\times \cL(_{r,r}M)$. Write $\cY$ for the stack quotient of $Y$ by the diagonal action of $\Sp(_{r,r}M)$. Let $\cA_{\cY}$ be the $\ZZ/2\ZZ$-graded purely of degree zero line bundle on $\cY$ with fibre $\det L_1\otimes\det L_2$ at $(L_1,L_2)$. Write $\tilde\cY$ for the gerb of square roots of $\cA_{\cY}$ over $\cY$. The map
$
\cL(_{r,r}M)\to Y
$
sending $L_1$ to $(_{0,r}M, L_1)\in Y$ yields a morphism of stacks 
$$
\rho: (\cL(_{r,r}M)/G(\cO/t^{2r}))\tilde{} \to \tilde\cY
$$
Write $S_{_{r,r}M}$ for the perverse sheaf on $\tilde\cY$ introduced in (Section~3.2, Definition~\ref{Def_S_M}). 
Set $\tau=\rho\comp\tilde\xi_r$.
It suffices to establish for any $r\ge 0$ a canonical isomorphism 
\begin{equation}
\label{iso_with_tau}
\tau^* S_{_{r,r}M}\otimes(\Qlb[1](\frac{1}{2}))^{\dimrel(\tau)}\,\iso\, \Aut
\end{equation}
over $_{rx}\Bunt_G$. 

 Remind that $Y_i\subset Y$ is the locally closed subscheme given by $\dim(L_1\cap L_2)=i$ for $(L_1,L_2)\in Y$. Let $\cY_i$ be the stack quotient of $Y_i$ by the diagonal action of $\Sp(_{r,r}M)$, set $\tilde \cY_i=\cY_i\times_{\cY} \tilde\cY$. Set 
$$
 _{rx,i}\Bunt_G={_{rx}\Bunt_G}\cap {_i\Bunt_G}\;\;\;\mbox{and} \;\;\; {_{rx,i}\Bun_G}={_{rx}\Bun_G}\cap {_i\Bun_G}
$$
For each $i$ the map $\tau$ fits into a cartesian square
$$
\begin{array}{ccc}
_{rx,i}\Bunt_G & \toup{\tau_i} & \tilde\cY_i\\
\downarrow && \downarrow\\
_{rx}\Bunt_G & \toup{\tau} & \tilde\cY
\end{array}
$$
Indeed, for $\cM\in {_{rx}\Bun_G}$ the space $\H^0(X,\cM)$ equals the intersection of $\cM/\cM(-rx)$ and $\H^0(X,\cM(rx))$ inside $\cM(rx)/\cM(-rx)$. 
By (\cite{L1}, Theorem~1), the $*$-restriction of $\Aut$ to $_i\Bunt_G\,\iso\,$ ${_i\Bun_G}\times B(\mu_2)$ identifies with 
$$
(\Qlb\boxtimes W)\otimes(\Qlb[1](\frac{1}{2}))^{\dim\Bun_G-i}
$$
Simialrly, by (\cite{L1}, Proposition~1 and 5), the $*$-restriction of $S_M$ to $\tilde\cY_i\,\iso\,\cY_i\times B(\mu_2)$ identifies with 
$$
(\Qlb\boxtimes W)\otimes(\Qlb[1](\frac{1}{2}))^{\dim\cY-i}
$$
Since the map $\tau_i$ is compatible with our trivializations of the corresponding gerbs, we get the isomorphism (\ref{iso_with_tau}) over $_{rx,i}\Bunt_G$ for each $i$. Since $\Aut$ is perverse, this also shows that the LHS of (\ref{iso_with_tau}) is placed in perverse degrees $\le 0$, and its $*$-restriction to $_{\le 2}\Bunt_G$ is placed in perverse degrees $<0$. 

 The map $\tau$ is not smooth, we overcome this difficulty as follows. 
Let us show that the LHS of (\ref{iso_with_tau}) is placed in perverse degrees $\ge 0$. Consider the stack $\cX$ classifying $(\cM,\cB)\in{_{rx}\Bunt_G}$ and a trivialization 
$$
\cM\mid_{\Spec\cO_x/t_x^{2r}}\,\iso\, M\mid_{\Spec\cO_x/t_x^{2r}}
$$ 
of the corresponding $G$-torsor. Let $\nu:\cX\to \tilde Y$ be the map
sending a point of $\cX$ to the triple $(\cM/\cM(-rx), H^0(X, \cM(rx)), \cB)$. Define $\cX_1$ and $\cX_3$ by the cartesian squares
$$
\begin{array}{ccc}
\cX_3 & \to & C_3\\
\downarrow\lefteqn{\scriptstyle\pi_{\cX_3}} && \downarrow\lefteqn{\scriptstyle\pi_C}\\
\cX_1 & \to & U\times_{\cL(_{r,r}M)} U\\
\downarrow && \downarrow\lefteqn{\scriptstyle \pi_Y}\\
\cX & \toup{\nu} & \tilde Y,
\end{array}
$$
Using (\ref{iso_from_AnnENS}), we get an isomorphism
$$
\mu^*\tau^*S_{_{r,r}M}\otimes(\Qlb[1]
(\frac{1}{2}))^{\dimrel(\mu)+\dimrel(\tau)}\,\iso\, (\pi_{\cX_3})_! \cE\otimes(\Qlb[1](\frac{1}{2}))^{\dim\cX_3}
$$
for some rank one local system $\cE$ on $\cX_3$. Here $\mu: \cX_1\to {_{rx}\Bunt_G}$ is the projection, it is smooth. Since $\pi_{\cX_3}$ is affine and $\cX_3$ is smooth, the LHS of (\ref{iso_with_tau}) is placed in perverse degrees $\ge 0$. 

 Thus, there exists an exact sequence of perverse sheaves $0\to K\to K_1\to\Aut\to 0$ on $_{rx}\Bunt_G$, where $K_1=\tau^* S_{_{r,r}M}\otimes(\Qlb[1](\frac{1}{2})) ^{\dimrel(\tau)}$, and $K$ is the extension by zero from $_{\le 2}\Bunt_G$. But we know already that $K_1$ and $\Aut$ are isomorphic in the Grothendieck group of $_{rx}\Bunt_G$. So, $K$ vanishes in this Grothendieck group, hence $K=0$. We are done.
\end{Prf}

\end{document}